 %
%

\input amstex
\documentstyle{amsppt}
\magnification=1200

 \vsize19.5cm \hsize13.5cm
 \TagsOnRight
 \pageno=1
\baselineskip=15.0pt
\parskip=3pt



\def\p{\partial}
\def\noo{\noindent}
\def\eps{\varepsilon}
\def\lam{\lambda}
\def\Om{\Omega}

\def\pom{{\p \Om}}
\def\bom{{\overline\Om}}
\def\R{\bold R}

\def\th{\theta}
\def\wtt{\widetilde}

\def\Ga{\Gamma}

\def\diam{\text{diam}}
\def\dist{\text{dist}}
\def\det{\text{det}}

\def\ol{\overline}

\def\D{\nabla}
\def\phi{\varphi}

\nologo
\NoRunningHeads

\topmatter

\title  {Boundary regularity for the Monge-Amp\`ere \\
               and affine maximal surface equations}\endtitle

\author{\bf Neil S. Trudinger \ \ \ \ Xu-Jia Wang}\endauthor

\affil{Centre for Mathematics and Its Applications\\
       The Australian National University\\
       }\endaffil

\address{Centre for Mathematics and Its Applications,
       Australian National University, \newline
       Canberra, ACT 0200,
       Australia}\endaddress

\email{Neil.Trudinger\@maths.anu.edu.au,\ \ \ \ \ \
X.J.Wang\@maths.anu.edu.au}\endemail

\thanks{Supported by Australian Research Council.
\newline
This paper was submitted for publication in September 2004.}
\endthanks

\abstract {In this paper, we prove global second derivative
estimates for solutions of the Dirichlet problem for the
Monge-Amp\`ere equation when the inhomogeneous term is only
assumed to be H\"older continuous. As a consequence of our
approach, we also establish the existence and uniqueness of
globally smooth solutions to the second boundary value problem for
the affine maximal surface equation and affine mean curvature
equation.}\endabstract

\endtopmatter

\document

\baselineskip=14.50pt
\parskip=2.5pt

\centerline {\bf \S 1. Introduction}

\vskip10pt

In a landmark paper [4], Caffarelli established interior $W^{2,
p}$ and $C^{2,\alpha}$ estimates for solutions of the
Monge-Amp\`ere equation
$$\det D^2 u  =f \tag 1.1$$
in a domain $\Om$ in Euclidean $n$-space, $\R^n$, under minimal
hypotheses on the function $f$. His approach in [3,4] pioneered
the use of affine invariance in obtaining estimates, which
hitherto depended on uniform ellipticity, [2,19], or stronger
hypotheses on the function $f$, [9, 13, 18]. If the function $f$
is only assumed positive and H\"older continuous in $\Om$, that is
$f\in C^\alpha (\Om)$ for some $\alpha\in (0, 1)$, then one has
interior estimates for convex solutions of (1.1) in $C^{2,
\alpha}(\Om)$ in terms of their strict convexity. When $f$ is
sufficiently smooth, such estimates go back to Calabi and
Pogorelov [9, 18]. The estimates are not genuine interior
estimates as assumptions on Dirichlet boundary data are needed to
control the strict convexity of solutions [4, 18].

Our first main theorem in this paper provides the corresponding
global estimate for solutions of the Dirichlet problem,
$$u=\phi\ \ \ \ \text{on}\ \ \pom. \tag 1.2$$

\proclaim{Theorem 1.1} Let $\Om$ be a uniformly convex domain in
$\R^n$, with boundary $\pom\in C^3$, $\phi\in C^3(\bom)$ and $f\in
C^\alpha(\bom)$, for some $\alpha\in (0, 1)$, satisfying $\inf
f>0$. Then any convex solution $u$ of the Dirichlet problem (1.1),
(1.2) satisfies the a priori estimate
$$\|u\|_{C^{2, \alpha}(\bom)}\le C,\tag 1.3$$
where $C$ is a constant depending on $n, \alpha$, $\inf f$,
$\|f\|_{C^\alpha(\bom)}$, $\pom$ and $\phi$.
\endproclaim

The notion of solution in Theorem 1.1, as in [4], may be
interpreted in the generalized sense of Aleksandrov [18], with
$u=\phi$ on $\pom$ meaning that $u\in C^0(\bom)$. However by
uniqueness, it is enough to assume at the outset that $u$ is
smooth. In [22], it is shown that the solution to the Dirichlet
problem, for constant $f>0$, may not be $C^2$ smooth or even in
$W^{2, p}(\Om)$ for large enough $p$, if either the boundary
$\pom$ or the boundary trace $\phi$ is only $C^{2,1}$. But the
solution is $C^2$ smooth up to the boundary (for sufficiently
smooth $f>0$) if both $\pom$ and $\phi$ are $C^3$ [22].
Consequently the conditions on $\pom$, $\phi$ and $f$ in Theorem
1.1 are optimal.

As an application of our method, we also derive global second
derivative estimates for the second boundary value problem of the
affine maximal surface equation and, more generally, its
inhomogeneous form which is the equation of prescribed affine mean
curvature. We may write this equation in the form
$$L[u] := U^{ij} D_{ij}w =f\ \ \ \  \text{in}\ \ \Om ,  \tag 1.4$$
where $[U^{ij}]$ is the cofactor matrix of the Hessian matrix $D^2
u$ of the convex function $u$ and
$$w   =[\det D^2 u]^{-(n+1)/(n+2)}. \tag 1.5$$
The second boundary value problem for (1.4), (as introduced in
[21]), is the Dirichlet problem for the system (1.4) (1.5), that
is to prescribe
$$u  =\phi\ \ \ \ \ \
   w  =\psi\ \ \ \  \text{on}\ \ \pom , \tag 1.6$$

We will prove

\proclaim{Theorem 1.2} Let $\Om$ be a uniformly convex domain in
$\R^n$, with $\pom\in C^{3,1}$, $\phi\in C^{3,1}(\bom)$, $\psi\in
C^{3,1}(\bom)$, $\inf_\Om \psi>0$ and $f\le 0, \in L^\infty(\Om)$.
Then there is a unique uniformly convex solution $u\in W^{4,
p}(\Om)$ (for all $1<p<\infty$) to the boundary value problem
(1.4)-(1.6). If furthermore $f\in C^\alpha(\bom)$, $\phi\in C^{4,
\alpha}(\bom)$,  $\psi\in C^{4, \alpha}(\bom)$, and $\pom\in C^{4,
\alpha}$  for some $\alpha\in (0, 1)$, then the solution $u\in
C^{4, \alpha}(\bom)$
\endproclaim

The condition $f\le 0$, (corresponding to non-negative prescribed
affine mean curvature [1, 17]), is only used to bound the solution
$u$. It can be relaxed to $f\le \delta$ for some $\delta>0$, but
it cannot be removed completely.

The affine mean curvature equation (1.4) is the Euler equation of
the functional
$$J[u]=A(u)-\int_\Om fu, \tag 1.7$$
where
$$A(u)=\int_\Om [\det D^2 u]^{1/(n+2)}. \tag 1.8$$
is the affine surface area functional. The natural or variational
boundary value problem for (1.4), (1.7) is to prescribe $u$ and
$\D u$ on $\pom$ and is treated in [21]. Regularity at the
boundary is a major open problem in this case.

Note that the operator $L$ in (1.4) possesses much stronger
invariance properties than its Monge-Amp\`ere counterpart (1.1) in
that $L$ is invariant under unimodular affine transformations in
$\R^{n+1}$, (of the dependent and independent variables).

Although the statement of Theorem 1.1 is reasonably succinct, its
proof is technically very complicated. For interior estimates one
may assume by affine transformation that a section of a convex
solution is of good shape, that is it lies between two concentric
balls whose radii ratio is controlled. This is not possible for
sections centered on the boundary and most of our proof is
directed towards showing that such sections are of good shape.
After that we may apply a similar perturbation argument to the
interior case [4]. To show sections at the boundary are of good
shape we employ a different type of perturbation which proceeds
through approximation and extension of the trace of the
inhomogeneous term $f$. The technical realization of this approach
constitutes the core of our proof. Theorem 1.1 may also be seen as
a companion result to the global regularity result of Caffarelli
[6] for the natural boundary value problem for the Monge-Amp\`ere
equation, that is the prescription of the image of the gradient of
the solution, but again the perturbation arguments are
substantially different.

The organization of the paper is as follows. In the next section,
we introduce our perturbation of the inhomogeneous term $f$ and
prove some preliminary second derivative estimates for the
approximating problems. We also show that the shape of a section
of a solution at the boundary can be controlled by its mixed
tangential-normal second derivatives. In Section 3, we establish a
partial control on the shape of sections, which yields $C^{1,
\alpha}$ estimates at the boundary for any $\alpha\in (0, 1)$,
(Theorem 3.1). In order to proceed further, we need a modulus of
continuity estimate for second derivatives for smooth data and
here it is convenient to employ a lemma from [8], which we
formulate in Section 4. In Section 5, we conclude our proof that
sections at the boundary are of good shape, thereby reducing the
proof of Theorem 1.1 to analogous perturbation considerations to
the interior case [4], which we supply in Section 6 (Theorem 6.1).
Finally in Section 7, we consider the application of our preceding
arguments to the affine maximal surface and affine mean curvature
equations, (1.4). In these cases, the global second derivative
estimates follow from a variant of the condition $f\in
C^\alpha(\bom)$ at the boundary, namely
$$|f(x)-f(y)| \le C|x-y|, \tag 1.9$$
for all $x\in\Om, y\in \pom$, which is satisfied by the function
$w$ in (1.5). The uniqueness part of Theorem 1.2 is proved
directly, (using an argument based on concavity), and the
existence part follows from our estimates and a degree argument.
The solvability of (1.4)-(1.6) without boundary regularity was
already proved in [21] where it was used to prove interior
regularity for the first boundary value problem for (1.4).

\vskip30pt

\centerline{\bf 2. Preliminary estimates}

\vskip10pt

Let $\Om$ be a uniformly convex domain in $\R^n$ with $C^3$
boundary, and $\phi$ be a $C^3$ smooth function on $\bom$. For
small positive constant $t>0$, we denote $\Om_t=\{x\in\Om\ |\
\dist(x, \pom)>t\}$ and $D_t=\Om-\bom_t$. For any point
$x\in\bom$, we will use $\xi$ to denote a unit tangential vector
of $\pom_\delta$ and $\gamma$ the unit outward normal of
$\pom_\delta$ at $x$, where $\delta=\dist(x, \pom)$.

Let $u$ be a solution of (1.1) (1.2). By constructing proper
sub-barriers we have the gradient estimate
$$\sup_{x\in\Om} |Du(x)|\le C.\tag 2.1 $$
We also have the second order tangential derivative estimates
$$C^{-1} \le u_{\xi\xi}(x)\le C\tag 2.2$$
for any $x\in\pom$. The upper bound in (2.2) follows directly from
(2.1) and the boundary condition (1.2). For the lower bound,  one
requires that $\phi$ is $C^3$ smooth, and $\pom$ is $C^3$ and
uniformly convex [22]. For (2.1) and (2.2) we only need that $f$
is a bounded positive function.

In the following we will assume that $f$ is positive and $f\in
C^\alpha(\bom)$ for some $\alpha\in (0, 1)$. Let $f_\tau$ be the
mollification of $f$ on $\pom$, namely $f_\tau=\eta_\tau * f$,
where $\eta$ is a mollifier on $\pom$. If $t>0$ is small, then for
any point $x\in D_t$, there is a unique point $\hat x\in\pom$ such
that $\dist(x, \pom)=|x-\hat x|$ and $\gamma=(\hat x-x)/|\hat
x-x|$. Let
$$f_t (x)=\cases
f(x)\ \ \ &\text{in}\ \  \Om_{2t},\\
f_\tau(\hat x)-C\tau^\alpha\ \  &\text{in}\ \ D_t,\\
\endcases \tag 2.3$$
where
$$\tau=t^{\eps_0}\ \ \ \ \eps_0=1/4n.$$
We define $f_t$ properly in the remaining part $\Om_t-\Om_{2t}$
such that, with a proper choice of the constant $C=C_{t,\tau}>0$,
$f_t\le f$ in $\Om$ and $f_t$ is H\"older continuous in $\bom$
with H\"older exponent $\alpha'=\eps_0\alpha$,
$$\align
|f_t-f|  \le C\tau^\alpha & =Ct^{\alpha'}
          \ \ \text{in}\ \ \Om, \\
\|f_t\|_{C^{\alpha'}(\bom)}&\le C\|f\|_{C^\alpha(\bom)}\\
\endalign$$
for some $C>0$ independent of $t$. From (2.3), $f_t$ is smooth in
$D_t$,
$$|D f_t| \le C\tau^{\alpha-1},\ \ \ \
 |D^2 f_t|\le C\tau^{\alpha-2},\ \ \
 \text{and}\ \ |\p_\gamma f_t| =0\ \ \ \text{in}\ \ D_t. \tag 2.4$$
Let $u_t$ be the solution of the Dirichlet problem,
$$\align
 \det D^2 u &= f_t\ \ \ \text{in}\ \ \Om,\tag 2.5\\
          u &= \phi\ \ \ \text{on}\ \ \pom.\\
\endalign $$
First we establish some a priori estimates for $u_t$ in $D_t$.
Note that by the local strict convexity [3] and the a priori
estimates for the Monge-Amp\`ere equation [18], $u$ is smooth in
$D_t$.

For any given boundary point, we may suppose it is the origin such
that $\Om\subset \{x_n>0\}$, and locally $\pom$ is given by
$$x_n=\rho(x')\tag 2.6$$
for some $C^3$ smooth, uniformly convex function $\rho$ satisfying
$\rho(0)=0, D \rho(0)=0$, where $x'=(x_1, \cdots, x_{n-1})$.  By
subtracting a linear function we also suppose that
$$u_t(0)=0,\ \ \ \ Du_t(0)=0. \tag 2.7$$
We make the linear transformation $T: x\to y$ such that
$$\align
y_i & = x_i/\sqrt t ,\ \ \ i=1, \cdots, n-1,\\
y_n & = x_n/ t , \tag 2.8\\
v & =u_t/t . \\
\endalign $$
Then $v$ satisfies  the equation
$$\det D^2 v=t f_t
 \ \ \ \text{in}\ \ T(\Om). \tag 2.9$$
Let $G=T(\Om)\cap\{y_n<1\}$. In $G$ we have $0\le v\le C$ since
$v$ is bounded on $\p G\cap\{y_n<1\}$. Observe that the boundary
of $G$ in $\{y_n<1\}$ is smooth and uniformly convex. Hence
$$|v_\gamma|  \le C
     \ \ \ \text{in}\ \p G\cap \{y_n<\frac 78\}.$$
From (2.2) we have
$$C^{-1}\le  v_{\xi\xi} \le C
     \ \ \ \text{on}\ \p G\cap \{y_n<\frac 78\}.$$
The mixed derivative estimate
$$|v_{\gamma\xi}| \le C
     \ \ \ \text{on}\ \p G\cap \{y_n<\frac 34\},$$
where $v_{\xi\gamma}=\sum \xi_i\gamma_j v_{y_iy_j}$, is found for
example in [8, 13]. For the mixed derivative estimate we need
$f_t\in C^{0, 1}$, with
$$|Df_t|\le C\tau^{\alpha-1}t^{1/2}\le C .$$
From (2.2) and equation (2.9) we have also
$$ v_{\gamma\gamma}\le C
      \ \ \ \text{on}\ \p G\cap \{y_n<\frac 34\}. $$
Next we derive an interior estimate for $v$.

\proclaim{Lemma 2.1} Let $v$ be as above. Then
$$|D^2 v|\le C(1+M)\ \ \ \text{in}\ G\cap \{y_n<\frac 12\},\tag 2.10$$
where $M=\sup_{\{y_n<7/8\}}|Dv|^2$, $C>0$ is independent of $M$.
\endproclaim

\noo{\it Proof}.  First we show $v_{ii}\le C$ for $i=1, \cdots,
n-1$. Let
$$w(y)=\rho^4 \eta(\frac 12 v_1^2)v_{11} , $$
where $v_1=v_{y_1}$, $v_{11}=v_{y_1y_1}$, and $\rho(y)=2-3y_n$ is
a cut-off function, $\eta(t) = (1-\frac tM )^{-1/8}$. If $w$
attains its maximum at a boundary point, by the above boundary
estimates we have $w\le C$. If $w$ attains its maximum at an
interior point $y_0$, by the linear transformation
$$\align
\wtt y_i & =y_i, \ \ \ i=2, \cdots, n , \\
\wtt y_1 & =y_1-\frac {v_{1i}(y_0)}{v_{11}(y_0)} y_i , \\
\endalign $$
which leaves $w$ unchanged, one may suppose $D^2 v(y_0)$ is
diagonal. Then at $y_0$ we have
$$\align
0 & =(\log w)_i = 4\frac {\rho_i}{\rho}
    + \frac {\eta_i}{\eta} + \frac {v_{11i}}{u_{11}} ,\tag 2.11\\
0 & \ge (\log w)_{ii} =
     4(\frac {\rho_{ii}}{\rho}-\frac {\rho_i^2}{\rho^2})
    + (\frac {\eta_{ii}}{\eta}-\frac {\eta_i^2}{\eta^2})
    + (\frac { v_{11ii}}{v_{11}}
                     -\frac {v_{11i}^2}{v_{11}^2} ) . \tag 2.12\\
\endalign $$
Inserting (2.11) into (2.12) in the form $\frac {\rho_i}{\rho}  =
-\frac 14(\frac {\eta_i}{\eta}
    + \frac {v_{11i}}{v_{11}})$ for $i=2, \cdots, n$ and
    $\frac {v_{11i}}{v_{11}}  = -(4\frac {\rho_i}{\rho}
    + \frac {\eta_i}{\eta})$ for $i=1$,
we obtain
$$\align
0 \ge & v^{ii}(\log w)_{ii} \\
  \ge & v^{ii} (\frac {\eta_{ii}}{\eta}- 3\frac {\eta_i^2}{\eta^2})
         -36 v^{11}\frac {\rho_1^2}{\rho^2}
    + v^{ii}\frac {v_{11ii}}{v_{11}}
 -\frac 32\sum_{i=2}^{n} v^{ii} \frac {v_{11i}^2}{v_{11}^2},\tag 2.13\\
\endalign $$
where $(v^{ij})$ is the inverse matrix of $(v_{ij})$.

It is easy to verify that
$$v^{ii}(\frac {\eta_{ii}}{\eta}- 3\frac {\eta_i^2}{\eta^2})
    \ge \frac {C}{M} v_{11} -\frac {C}{M}, $$
where $C>0$ is independent of $M$. Differentiating the equation
$$\log \det D^2 v=\log(tf_t)$$
twice with respect to $y_1$, and observing that $|\p_1f_t|\le
C\tau^{\alpha-1}t^{1/2}\le C$ and $|\p_1^2 f_t|\le
C\tau^{\alpha-2}t\le C$ after the transformation (2.8), we see the
last two terms in (2.13) satisfy
$$v^{ii}\frac {v_{11ii}}{v_{11}}
 -\frac 32\sum_{i=2}^{n} v^{ii} \frac {v_{11i}^2}{v_{11}^2}
 \ge -\frac {1}{v_{11}}(\log f_t)_{11}\ge -C. $$
We obtain
$$\rho^4 v_{11} \le C(1+M) .$$
Hence $v_{ii}\le C$ for $i=1, \cdots, n-1$ in $G\cap \{y_n<\frac
12\}$.

Next we show that $v_{nn}\le C$.  Let $w(y)=\rho^4 \eta(\frac 12
v_n^2)v_{nn}$ with the same $\rho$ and $\eta$ as above. If $w$
attains its maximum at a boundary point, we have $v_{nn}\le C$ by
the boundary estimates. Suppose $w$ attains its maximum at an
interior point $y_0$. As above we introduce a linear
transformation
$$\align
\wtt y_i & =y_i, \ \ \ i=1, \cdots, n-1,\\
\wtt y_n & =y_n-\frac {v_{in}(y_0)}{v_{nn}(y_0)} y_i,\\
\endalign $$
which leaves $w$ unchanged. Then
$$w(y)=(2-\alpha_iy_i)^4 \eta(\frac 12 v_n^2)v_{nn} $$
and $D^2 v(y_0)$ is diagonal. By the estimates for $v_{ii}$, $i=1,
\cdots, n-1$, the constants $\alpha_i$ are uniformly bounded.
Therefore the above argument applies. $\square$

Scaling back to the coordinates $x$, we therefore obtain
$$\align
\p_\xi^2 u_t (x) & \le C
                \ \ \ \ \ \ \text{in}\ \ D_{t/2}, \tag 2.14a\\
|\p_\xi\p_\gamma u_t (x)| & \le C/\sqrt t
                      \ \ \ \text{in}\ \ D_{t/2}, \tag 2.14b\\
\p_\gamma^2 u_t (x) & \le C/t
                \ \ \ \ \ \ \text{in}\ \ D_{t/2}, \tag 2.14c\\
\endalign $$
where $C$ is independent of $t$, $\xi$ is any unit tangential
vector to $\pom_\delta$ and $\gamma$ is the unit normal to
$\pom_\delta$ ($\delta=\dist(x, \pom)$), and $\p_\xi\p_\gamma
u=\sum\xi_i\gamma_j u_{x_ix_j}$.

The proof of Lemma 2.1 is essentially due to Pogorelov [18]. Here
we used a different auxiliary function, from which we obtain a
linear dependence of $\sup|D^2 v|$ on $M$, which will be used in
the next section. The linear dependence can also be derived from
Pogorelov's estimate by proper coordinate changes. Taking
$\rho=-u$ in the auxiliary function $w$, we have the following
estimate.

\proclaim{Corollary 2.1} Let $u$ be a convex solution of $\det D^2
u= f$ in $\Om$. Suppose $\inf_\Om u=-1$, and either $u=0$ or $|D^2
u|\le C_0(1+M)$ on $\pom$. Then
$$|D^2 u|(x)\le C(1+M),\ \ \forall\ x\in \{u<-\frac 12\},\tag 2.15$$
where $M=\sup_{\{u<0\}} |Du|^2$, and $C$ is independent of $M$.
\endproclaim

Next we derive some estimates on the level sets of the solution
$u$ to (1.1) (1.2). Denote
$$\align
S^0_{h, u}(y) &
  = \{x\in\bom\ {|}\  u(x)<u(y)+Du(y)(x-y)+h\} ,\\
S_{h, u}(y) &
  = \{x\in\bom\ {|}\  u(x)=u(y)+Du(y)(x-y)+h\} .\\
\endalign $$
We will write $S_{h, u}=S_{h, u}(y)$ and $S^0_{h, u}=S^0_{h,
u}(y)$ if no confusion arises. The set $S^0_{h, u}(y)$ is the
section of $u$ at center $y$ and height $h$ [4].

\proclaim{Lemma 2.2} There exist positive constants $C_2>C_1$
independent of $h$ such that
$$C_1h^{n/2}\le |S^0_{h, u}(y)|\le C_2 h^{n/2} \tag 2.16$$
for any $y\in\pom$, where $|\Cal K|$ denotes the Lebesgue measure
of a set $\Cal K$.
\endproclaim

\noo{\it Proof}. It is known that for any bounded convex set $\Cal
K\subset \R^n$,   there is a unique ellipsoid $E$ containing $\Cal
K$ which achieves the minimum volume among all ellipsoids
containing $\Cal K$ [3]. $E$ is called the minimum ellipsoid of
$\Cal K$. It satisfies $\frac 1n (E-x_0)\subset \Cal K-x_0\subset
E-x_0$, where $x_0$ is the center of $E$.

Suppose the origin is a boundary point of $\Om$, $\Om\subset
\{x_n>0\}$, and locally $\pom$ is given by (2.6). By subtracting a
linear function we also suppose $u$ satisfies (2.7). Let $E$ be
the minimum ellipsoid of $S^0_{h, u}(0)$. Let $v$ be the solution
to $\det D^2 u=\inf_\Om f_t$ in $S^0_{h, u}$, $v=h$ on $\p S^0_{h,
u}$. If $|E|>Ch^{n/2}$ for some large $C>1$, we have $\inf v<0$.
By the comparison principle, we obtain $\inf u\le \inf v<0$, which
is a contradiction to (2.7). Hence the second inequality of (2.16)
holds.

Next we prove the first inequality. Denote
$$\align
a_h & =\sup\{|x'|\ {|}\ x\in S_{h, u}(0)\}, \tag 2.17\\
b_h & =\sup\{x_n\ {|}\ x\in S_{h, u}(0)\}. \tag 2.18\\
\endalign $$
If the first inequality is not true, we have $|S^0_{h,
u}|=o(h^{n/2})$ for a sequence $h\to 0$. By (2.2), we have
$S^0_{h, u}\supset \{x\in\pom\ | \ |x|<C h^{1/2}\}$ for some
$C>0$. Hence $b_h=o(h^{1/2})$. By (2.2) we also have $u(x)\ge
C_0|x|^2$ for $x\in \pom$. Hence if $a_h\le Ch^{1/2}$ for some
$C>0$, the function
$$v= \delta_0 (|x'|^2+ (\frac {h^{1/2}}{b_h}x_n)^2)+ \eps x_n$$
for some small $\delta_0>0$, is a sub-solution to the equation
$\det D^2 u=f$ in $S^0_{h, u}$ satisfying $v\le u$ on $\p S^0_{h,
u}$, where $\eps>0$ can be arbitrarily small. It follows by the
comparison principle that $v_n(0)\le u_n(0)=0$, which contradicts
with $v_n(0)=\eps>0$.

Hence we have $a_h/h^{1/2}\to\infty$ as $h\to 0$. Let
$x_0=(x_{0,1}, 0, \cdots, 0, x_{0, n})$ (after a rotation of the
coordinates $x'$) be the center of $E$, where $E$ is the minimum
ellipsoid of $S^0_{h, u}$.  Make the linear transformation
$$y_1 = x_1-(x_{0,1}/x_{0,n})x_n,\ \ \ \
 y_i = x_i\ \ \ \ i=2, \cdots, n$$
such that the center of $E$ is moved to the $x_n$-axis. Let
$E'=\{\sum_{i=1}^{n-1} (x_i/a_i)^2 <1\}$ be the projection of $E$
on $\{x_n=0\}$. Since the origin $0\in S^0_{h, u}$ and the center
of $E$ is located on the $x_n$-axis, one easily verifies that
$a_1\cdots a_n\le C|S^0_{h, u}|=o(h^{n/2})$, where $a_n=x_{0, n}$.
Note that $x_{0,1}\le a_h$ and $x_{0,n}\le b_h\le 2n x_{0,n}$. By
the uniform convexity of $\pom$, we have
$$\frac {x_{0, n}}{x_{0, 1}}\ge C\frac {b_h}{a_h}\ge Ca_h>>h^{1/2}.$$
Hence after the above transformation, the boundary part $\pom\cap
S^0_{h, u}$ is still uniformly convex. Hence as above the function
$v=\delta_0\sum_{i=1}^n (\frac {h^{1/2}}{a_i}y_i)^2+\eps y_n$ is a
sub-solution, and we also reach a contradiction. $\square$

Next we show that the shape of the level set $S_{h, u}$ can be
controlled by the mixed derivatives $u_{\xi\gamma}$ on $\pom$.

\proclaim{Lemma 2.3} Let $u$ be the solution of (1.1) (1.2).
Suppose as above that $\pom$ is given by (2.6) and $u$ satisfies
(2.7). If
$$|\p_{\xi\gamma} u (x)|\le K\ \ \ \text{on}\ \ \pom\tag 2.19$$
for some $K\ge 1$, then we have
$$\align
a_h & \le  CK h^{1/2},  \tag 2.20\\
b_h & \ge   C h^{1/2}/K \tag 2.21\\
\endalign $$
for some $C>0$ independent of $u$, $K$ and $h$.
\endproclaim

\noo{\it Proof}. We need only to prove (2.20) and (2.21) for small
$h>0$. Suppose the supremum $a_h$ is attained at $x_h =(a_h, 0,
\cdots, 0, c_h)\in S_{h, u}(0)$. Let $\ell=S_{h, u} \cap
\{x_2=\cdots=x_{n-1}=0\}$. Then $\ell\subset \bom$ and it has an
endpoint $\hat x=(\hat x_1,0, \cdots, 0,\hat x_n) \in \pom$ with
$\hat x_1>0$ such that $u(\hat x)=h$. If $a_h=\hat x_1$, by (2.2)
we have $\hat x_1\le Ch^{1/2}$, and by the upper bound in (2.16)
we have $b_h\ge Ch^{1/2}$. Hence (2.20) and (2.21) hold.

If $a_h>\hat x_1$, let $\xi=(\xi_1, 0, \cdots, 0, \xi_n)$ be the
unit tangential vector of $\pom$ at $\hat x$ in the
$x_1x_n$-plane, and $\zeta=(\zeta_1, 0, \cdots, 0, \zeta_n)$ be
the unit tangential vector of the curve $\ell$ at $\hat x$. Then
all $\xi_1, \xi_n, \zeta_1, $ and $\zeta_n>0$. Let $\th_1$ denote
the angle between $\xi$ and $\zeta$ at $\hat x$, and $\th_2$ the
angle between $\xi$ and the $x_1$-axis. By (2.2) and (2.19),
$$\align
|\p_\gamma u (\hat x)| & \le CK |\hat x|, \\
|\p_\xi u (\hat x)| & \ge C|\hat x|. \\
\endalign $$
Hence
$$\frac CK\le \th_1 < \pi-\frac CK.  \tag 2.22$$
But since all $\xi_1, \xi_n, \zeta_1, $ and $\zeta_n>0$, we have
$\th_1+\th_2< \frac \pi 2$. Note that by (2.2) and (2.16), $a_h\ge
Ch^{1/2}$ and $b_h\le Ch^{1/2}$. We obtain
$$\align
a_h & \le \hat x_1+ b_h/\text{tg}\, (\th_1+\th_2)
                                      \le CK h^{1/2} , \tag 2.23\\
b_h & \ge a_h\text{tg}\, (\th_1+\th_2) \ge  C h^{1/2}/K .\\
\endalign $$
Lemma 2.3 is proved. $\square$

Lemma 2.3 shows that the shape of the sections $S^0_{h, u}(y)$ at
boundary points $y$ can be controlled by the mixed second order
derivatives of $u$. If $S^0_{h, u}$ has a {\it good shape} for
small $h>0$, namely if the inscribed radius $r$ is comparable to
the circumscribed radius $R$,
$$R\le C_0r\tag 2.24$$
for some constant $C_0$ under control, the perturbation argument
[4] applies and one infers that $|D^2 u(0)|$ is bounded. See
Section 6. It follows that $u\in C^{2, \alpha}(\bom)$ by [2,19].
Estimation of the mixed second order derivatives on the boundary
will be the key issue in the rest of the paper.

\vskip30pt

\centerline{\bf 3. Mixed derivative estimates at the boundary }

\vskip10pt

For $t>0$ small let $u_t$ be a solution of (2.5) and assume (2.6)
(2.7) hold. As in Section 2 we use $\xi$ and $\gamma$ to denote
tangential (parallel to $\pom$) and normal (vertical to $\pom$)
vectors.

\proclaim{Lemma 3.1}  Suppose
$$|\p_\xi\p_\gamma u_t|\le K\ \ \ \text{on}\ \ \pom \tag 3.1$$
for some $1\le K\le Ct^{-1/2}$. Then we have
$$\align
\p_i^2 u_t\le C \ \ \ &\text{in}\
 D_t\cap\{x_n<t/8\} ,\ \ \ i=1, \cdots, n-1, \tag 3.2a\\
|\p_i\p_n u_t|\le C K
      \ \ \ &\text{in}\ D_t\cap\{x_n<t/8\} ,\tag 3.2b\\
\p_n^2 u_t\le CK^2 \ \ \
    &\text{in}\ D_t\cap\{x_n<t/8\} , \tag 3.2c\\
\endalign $$
where $C>0$ is a constant independent of $K$ and $t$.
\endproclaim

\noo{\it Proof.}\ By (2.14c), estimate (3.2a) is equivalent  to
(2.14a). The estimate (3.2b) follows from (3.2a) and (3.2c) by the
convexity of $u_t$. By (2.2), (3.1),  and equation (2.5), we
obtain (3.2c) on the boundary $\pom$. By (2.15), the interior part
of (3.2c) will follow if we have an appropriate gradient estimate
for $u_t$ in the set $S^0_{h, u_t}(0)$.

Let $h>0$ be the largest constant such that $S^0_{h,
u_t}(0)\subset D_{t/2}$ and $u_t$ satisfies (2.14) in $\{u_t<h\}$.
By the Lipschitz continuity of $u$, we have $h\le Ct$. Let $v(y) =
u_t(x)/h$, where $y = x/\sqrt h$. Then $v$ satisfies the equation
$$\det D^2 v=f_t\ \ \ \text{in}
                  \ \ \wtt\Om=\{x/\sqrt h\ |\  x\in\Om\}. \tag 3.3$$
By (2.16) we have
$$C_1\le |\{v<1\}|\le C_2.\tag 3.4$$
We claim
$$|\p_n v(y)|\le CK
         \ \ \ \forall\ \  y\in \{v<\frac 12\}. \tag 3.5$$
If (3.5) holds, by Corollary 2.1 (with the auxiliary function
$w(y)=(\frac 12-v)^4 \eta(\frac 12 v_n^2)v_{nn}$ in the proof of
Lemma 2.1), we obtain
$$\p_{y_n}^2 v\le C K^2\ \ \ \text{in}\ \ \{v<1/4\}.$$
In the above estimate we have used
$$\p_{y_n}^2 \log f_t(y)
   =h\,\p_{x_n}^2\log f_t(x)\le C \ \ \ \text{in}\ \{x_n<t\}$$
by our definition of $f_t$ in (2.3). Changing back to the
$x$-coordinates we obtain (3.2c).

By convexity it suffices to prove (3.5) for $ y\in \p\{v<\frac
12\}$. Let $\ol a_h=h^{-1/2}a_h$, where $a_h$ is defined in
(2.17). If $\ol a_h\le C$, by (2.16),  the set $\{v<1\}$ has a
good shape. By (2.1) and (2.2), the gradient estimate in
$\{v<\frac 12\}$ is obvious.

If $\ol a_h>>1$ ($\ol a_h\le CK$ by (2.20)), we divide
$\p\{v<\frac 12\}$ into two parts. Let $\p_1\{v<\frac 12\}$ denote
the set $y\in \{v=\frac 12\}\cap\wtt\Om$ such that the outer
normal line of $\{v<\frac 12\}$ at $y$ intersects with
$\{v=1\}=\{y\in\wtt\Om\ |\ v(y)=1\}$, and $\p_2\{v<\frac 12\}$
denote the rest of $\p\{v<\frac 12\}$, which consists of the
boundary part $\{v<\frac 12\}\cap\p\wtt\Om$ and the points $y\in
\{v=\frac 12\}$ at which the outer normal line of $\{v<\frac 12\}$
intersects with a boundary point in $\{v<1\}\cap\p\wtt\Om$.

Observe that for any $y\in \{v<1\}\cap\p\wtt\Om$, (3.5) holds by
(3.1) since $Dv(0)=0$. By convexity we obtain (3.5) on the part
$\p_2\{v<\frac 12\}$.

To verify (3.5) on $\p_1\{v<\frac 12\}$, we will construct
appropriate sub-barriers to show that the distance from $\{v=1\}$
to $\{v<\frac 12\}$ is greater than $C/K$. Then by the convexity
of $v$ we have $|Dv|<CK$ on $\p_1\{v<\frac 12\}$.

Our sub-barrier will be a function defined on a cylinder
$U=E\times(-a_n, a_n)\subset \R^n$, where $E=\sum_{i=1}^{n-1}
{x_i^2}/{a_i^2}<1$ is an ellipsoid in $\R^{n-1}$. Suppose
$a_1\cdots a_n=1$.  Let $w$ be the convex solution to the
Monge-Amp\`ere equation $\det D^2 w=1$ in $U$ satisfying $w=0$ on
$\p U$.

By making the linear transformation $\wtt y_i=y_i/a_i$ for $i=1,
\cdots, n$ such that $U=\{|{\wtt y}'|<1\}\times (-1, 1)$, where
${\wtt y}'=(\wtt y_1, \cdots, \wtt y_{n-1})$, we have the estimate
$C_1\le -\inf_U w\le C_2$ for two constants $C_2>C_1>0$ depending
only on $n$. By constructing proper sub-barriers [4], we see that
$w$ is H\"older continuous in $\wtt y$. Hence for any $C_0>0$, by
the convexity of $w$ we have the gradient estimate $C_1<|D_{\wtt
y}w|<C_2$ on $\{w<-C_0\}$, for different $C_2>C_1>0$ depending
only on $n$ and $C_0$. Changing back to the variable $y$, we
obtain
$$C_1a_n^{-1}\le |D_{y_n}w|\le C_2 a_n^{-1}\tag 3.6$$
at any point $y\in \{w=-C_0\}$ such that $y'\in \frac {1}{2} E$.
If $a:=a_1\cdots a_n\ne 1$, then by a dilation one sees that (3.6)
holds with $a_n$ replaced by $a_n/a$.

In order to use (3.6) to verify (3.5) on the part $\p_1\{v<\frac
12\}$, we first show that
$$\inf_{|\nu|=1}\sup_{y, z\in\{v<1\}}\nu\cdot (y-z)\ge C/K,
       \tag 3.7$$
namely the in-radius of the convex set $\{v<1\}$ is greater than
$C/K$, where $\nu\cdot y$ denotes the inner product in $\R^n$. To
prove (3.7) we first observe that by (2.2),
$$B_{r_1}(0)\cap\p\wtt\Om  \subset \{v<1\}\cap\p\wtt\Om
                 \subset B_{r_2}(0)\cap\p\wtt\Om \tag 3.8$$
for some $r_1, r_2>0$ independent of $t$. Let $\wtt y=(0, \cdots,
0, \wtt y_n)$ be a point on the positive $x_n$-axis such that
$v(\wtt y)=1$. To prove (3.7), it suffices to show that
$$\wtt y_n\ge C/K.\tag 3.9$$
Let $\ol y=(\ol a,0,\cdots, 0,\ol c)\in \p\wtt\Om$ be an arbitrary
point such that $v(\ol y)=1$. Then similarly to (2.22), the angle
at $\ol y$ of the triangle with vertices $\ol y, \wtt y$ and the
origin is larger than $C/K$. Hence $\wtt y_n\ge Cr_1/K\ge C/K$.
Hence (3.9) holds.

Now for any given point $\hat y\in\{v=1\}\cap\p\wtt\Om$, let $P$
denote the tangent plane of $\{v=1\}$ at $\hat y$. Choose a new
coordinate system $z$ such that $\hat y$ is the origin,
$P=\{z_n=0\}$ and the inner normal of $\{v<1\}$ is the positive
$z_n$-axis. Let $S'$ denote the projection $\{v<1\}$ on $P$. By
(3.4) and  (3.7) we have the volume estimate
$$|S'|\le CK. \tag 3.10$$
Let $E\subset P$ be the minimum ellipsoid of $S'$ with center
$z_0$, and $E_0\subset P$ be the translation of $E$ such that its
center is located at the origin $z=0$ (the point $\hat y$). Then
we have $S'\subset E\subset 4nE_0$. The latter inclusion is true
when $E$ is a ball and it is also invariant under linear
transformations.

Let $U=\beta E_0\times (0, 2/K)$ and $U_{1/2}=\beta E_0\times (0,
1/K)$. Let $w$ be the solution of $\det D^2 w=\sup_\Om f_t$ in $U$
such that $w=1$ on $\p U$. We may choose the constant $\beta\ge
8n$ such that $2E\subset \beta E_0$ and $\inf_U w \le -1$ (Note
that since $|U|=2\beta^{n-1}|E_0|/K$,  $\beta$ can be very large
if $|E_0|<<K$). Then by convexity we see that $w\le 0\le v$ on
$\{z_n=1/K\}\cap\{v<1\}$.

To verify that $w<v$ on $\p\wtt\Om\cap \{v<1\}$, we observe that
either the distance from the plane $P=\{z_n=0\}$ to the set
$\{v<1\}\cap\p\wtt\Om$ is larger than $C/K$, or the angle $\th_1$
between the plane $P$ and the plane $\{y_n=0\}$ satisfies (2.22).
In the former case, by (3.6) (with $a_n=1/K$) we have $w\le v$ on
$\p\wtt\Om\cap U_{1/2}$ if $\beta$ is chosen large, independent of
$K$. In the latter case, noting that the boundary part
$\p\wtt\Om\cap\{v<1\}$ is very flat and that $|\p_\xi v|\le C$,
where $\xi$ is tangential to $\p\wtt\Om$, by (3.6) we also have
$w\le v$ on $\p\wtt\Om\cap U_{1/2}$. Therefore in both cases we
have $w\le v$ on the boundary of the set $\{v<1\}\cap U_{1/2}$.

By the comparison principle, it follows that $w\le v$ in
$\{v<1\}\cap U_{1/2}$. By the gradient estimate (3.6) for $w$, it
follows that the distance from $\{v<\frac 12\}$ to  $\{v=1\}$ is
greater than $C/K$. This completes the proof. $\square$

\proclaim{Lemma 3.2}  Suppose $|D^2 u_t|\le K^2$ in $D_{t/8}$.
Then
$$|D^2 u_t|\le CK^2  \ \ \text{in}\ D_{2t}. \tag 3.11 $$
where $C>0$ is a constant independent of $K$ and $t$.
\endproclaim

\noo{\it Proof}. Fix a point $x_0\in D_{2t}-D_{t/8}$. For any
small $h>0$, there exists a linear function $x_{n+1}=a\cdot x+b$
such that $a\cdot x_0+b=u(x_0)+h$ and $x_0$ is the center of the
minimum ellipsoid $E$ of the section $\hat S_h :=\{x\in\Om\ | \
u(x)<a\cdot x+b\}$ [5], where $a$ and $b$ depend on $h$. Let $h$
be the largest constant such that $\hat S_{h-\eps}\subset
\subset\Om$ for any $\eps>0$.

Make a linear transformation $y=Tx$ such that $T(E)$ is a unit
ball. Let $v=|T|^{2/n}(u-a\cdot x-b)$. Then $v$ satisfies the
equation $\det D^2 v=f_t(T^{-1}(y))$ in $T(\hat S_h)$ and $v=0$ on
the boundary $\p T(\hat S_h)$. We have $C_1\le -\inf v\le C_2$ for
two constants $C_2>C_1>0$ depending only on $n$, the upper and
lower bounds of $f_t$. Let us assume simply that $\inf v=-1$.

Since $f_t$ is H\"older continuous with exponent
$\alpha'=\eps_0\alpha$, both before and after the transformation,
by the Schauder type estimate [4], we have $u\in C^{2,
\alpha'}(T(\hat S_h))$. That is for any $\delta>0$, there exist
$C_2>C_1>0$ depending on $n, \delta$, $\alpha'\in (0, 1)$, the
upper and lower bounds of $f_t$, and
$\|f_t\|_{C^{\alpha'}(\bom)}$, but independent of $h$, such that
$$C_1I\le \{D^2_yv(y)\}\le C_2 I  \tag 3.12$$
for any $y\in \{v<-\delta\}$, where $I$ is the unit matrix. Note
that (3.12) implies that the largest eigenvalue of $\{D^2_yv\}$ is
controlled by the smallest one.

Let $\delta=1/64$. Since $\inf v=-1$, by convexity, $v(y_0)\le
-\frac 12$,  where $y_0=T(x_0)$. Since $\dist(x_0, \pom)\le 2t$,
by convexity, there exists a point $x^*\in D_{t/8}$ such that
$v(y^*)\le -1/64$, where $y^*=T(x^*)$. From (3.12) we have
$$|D^2_yv(y_0)|\le C|D^2_y v(y^*)|.     $$
Changing back to the $x$-variables, we obtain (3.11). $\square$

\vskip10pt

The next lemma is simple but is important for our proof.

\proclaim {Lemma 3.3} Suppose
$$|D^2 u_t|\le C_0t^{\beta-1}\ \ \ \text{in}\ \ D_{2t},\tag 3.13$$
where $\beta\in [0, 1]$ is a constant. Then in $D_{t/2}$, we have
$$|u_t-u|(x)\le Ct^{\beta+\alpha'}\dist(x, \pom) , \tag 3.14$$
where $\alpha'=\eps_0\alpha$, $C$ is independent of $t$.
\endproclaim

\noo{\it Proof}. By our construction we have $f_t\le f$ in $\Om$.
Hence $u_t\ge u$ in $\Om$. Let
$$z=\cases
- 4t^{\beta+\alpha'} d_x + t^{\beta+\alpha'-1} d_x^2
                                \ \ \  &\text{if}\ d_x<2t,\\
- 4t^{\beta+\alpha'+1}  &\text{if}\ d_x\ge 2t,\\
\endcases \tag 3.15$$
where $d_x=\dist(x, \pom)$. For any point $x\in D_{2t}$, choose
the coordinates properly such that $D^2 z$ is diagonal with
$z_{11}\le \cdots\le z_{nn}$. Then
$$\det D^2 (u_t + C'z)
 \ge \det D^2 u_t +C'(\det \wtt D^2 u_t) z_{nn}, $$
where $\wtt D^2 u=(u_{ij})_{i, j=1}^{n-1}$. From (3.15) we have
$z_{nn}\ge C t^{\beta+\alpha'-1}$. By (3.13) we have  $\det \wtt
D^2 u_t\ge Ct^{1-\beta}$. Hence
$$\det D^2 (u_t + C'z)
 \ge f_t +C't^{\alpha'}\ge f $$
if $C'$ is chosen large. By the comparison principle, we obtain
(3.14). $\square$

In Lemma 3.2 we assume that $f\in C^\alpha(\bom)$ for some
$\alpha\in (0, 1)$. This condition is not satisfied in the proof
of Theorem 1.2. We will need the following alternative of Lemma
3.3 in this case.

\proclaim{Lemma 3.3$'$} Suppose $f$ satisfies
$$|f(x)-f(y)|\le C|x-y|
     \ \ \ \forall\ \ x\in\Om, y\in\pom.\tag 3.16$$
Then we have
$$|u-u_t|(x) \le Ct^{1+1/n} \dist(x, \pom) \tag 3.17$$
for some constant $C>0$ independent of $t$.
\endproclaim

\noo{\it Proof.}\  Let
$$z=\cases
- 4t^{1+1/n} d_x + t^{1/n} d_x^2 \ \ \  &\text{if}\ d_x<2t,\\
- 4t^{2+1/n}  &\text{if}\ d_x\ge 2t.\\
\endcases\tag 3.18$$
We have
$$\det D^2 z\ge Ct^n\ \ \ \text{in}\ \ D_{2t}$$
for some $C>0$. Under assumption (3.16), we have $|f_t-f|\le Ct$.
Hence
$$\det D^2 (u_t + Cz)
 \ge \det D^2 u_t +C(\det D^2 u_t)^{(n-1)/n} (\det D^2 z)^{1/n}
 \ge f\ \ \text{in}\ \ \Om. $$
Similarly we have $\det D^2 (u + Cz) \ge \det D^2 u_t$ in $\Om$.
It follows
$$|u-u_t|(x)\le C|z(x)|.$$
Hence (3.17) holds. $\square$

\vskip10pt

Let $\th=\alpha/16n$ if $f\in C^\alpha$, or $\th=1/16n$ if $f$
satisfies (3.16), and $t'=t^{1+\th}$. Let $u_{t'}$ be the
corresponding solution of (2.5). By our construction of $f_t$, we
may assume that $f_{t'}\ge f_t$ so that $u_{t'}\le u_t$. Obviously
Lemma 3.3 holds with $u$ replaced by $u_{t'}$.

\proclaim {Lemma 3.4} Suppose $u_t$ satisfies (3.1). Then
$$|\p_\xi\p_\gamma u_{t'}|\le CK
               \ \ \ \text{on}\ \ \pom, \tag 3.19$$
where $C$ is independent of $K$ and $t$.
\endproclaim

\noo{\it Proof.}\ Suppose the origin is a boundary point and (2.6)
(2.7) hold. For any $(x', s)\in\Om$, where $s=t'/8$, we have
$$\align
\p_i u_t (x', s) & =\p_i u_t(x', \rho(x'))
            +\p_n\p_i u_t(x', s_1)(s-\rho(x')),\ \ i<n \tag 3.20 \\
\p_i \phi (x', s) & =\p_i \phi (x', \rho(x'))
            +\p_n\p_i \phi (x', s_2)(s-\rho(x')) , \\
\endalign $$
for some $s_1, s_2\in (\rho(x'), s)$. Since $Du_t(0)=0$, by (3.1)
we have $|\p_\gamma u_t(x', \rho(x'))|\le CK|x'|$. Hence
$$|\p_n u_t(x', \rho(x'))|\le CK|x'|. $$
Since $\p_\xi (u_t-\phi)=0$, we obtain
$$|\p_i(u_t-\phi)|(x', \rho(x'))\le C |x'|\, |\p_n (u_t-\phi)|
      \le CK|x'|^2\le CKs.$$
By (3.2b) and (3.20) it follows that
$$ |\p_i (u_t-\phi) (x', s)|\le CKs . \tag 3.21$$

Let $\beta\in [0, 1]$ such that $K=t^{(\beta-1)/2}$ (by (2.14c) we
may assume $\beta\le 1$). Then by (3.1) and Lemmas 3.1 and 3.2,
$|D^2 u_t|\le Ct^{\beta-1}$ in $D_{2t}$. Hence by Lemma 3.3,
$$|u_t-u_{t'}| \le C t^{\beta+\alpha'} s
                   \ \ \ \text{on}\ \ \Om\cap \{x_n=s\} . $$
By (3.2a),
$$\p_i^2 u_t  \le C \ \ \ \text{and}\ \ \
\p_i^2 u_{t'}  \le C
   \ \ \ \text{on}\ \ \Om\cap \{x_n=s\} .$$
Hence
$$|\p_i (u_t-u_{t'})|
 \le C\sup_{\{x_n=s\} }|u_t-u_{t'}|^{1/2}
  \le C (t^{\beta+\alpha'} s)^{1/2}
    \ \ \ \ \ \text{on}\ \ \Om\cap \{x_n=s\} .$$
Recall that $s=t'/8=t^{(1+\th)}/8$. We obtain
$$\align
|\p_i (u_t-u_{t'})|
 &\le Ct^{(\beta+\alpha'-1-\th)/2} s\\
  &\le Ct^{(\beta-1)/2} s=CKs
        \ \ \ \ \ \text{on}\ \ \Om\cap \{x_n=s\}.  \tag 3.22\\
        \endalign $$
From (3.21) and (3.22) we thus obtain
$$|\p_i(\phi - u_{t'})|\le CKs
    \ \ \ \text{on}\ \ \{x_n=s\}. \tag 3.23$$
Next we estimate $\p_nu_{t'}$ on $\{x_n=s\}$. First we consider
the point $(0, s)$. By convexity and (3.14) we have
$$\align
\p_nu_{t'}(0, s)
  &\le \frac 1s[u_{t'}(0, 2s)-u_{t'}(0, s)]\\
  &\le \frac 1s[u_t(0, 2s)-u_t (0, s)]+Cs^{\beta+\alpha'}\\
   &\le \p_n u_t(0, 2s)+Cs^{\beta+\alpha'} . \\
   \endalign $$
By Lemma 3.1, $\p_n^2 u_t\le CK^2$. Hence $\p_n u_t(0, 2s)\le
\p_nu_t(0)+CK^2 s=CK^2s$. Note that $Ks^{1/2}< Kt^{(1+\th)/2}$
$\le t^{\th/2}$. Hence we obtain
$$\align
\p_nu_{t'}(0, s)
   & \le CK^2s +Cs^{\beta+\alpha'}\\
   &\le Ct^{\th/2} K s^{1/2}+ Cs^{\beta+\alpha'}
     \le CKs^{1/2} .\\
\endalign $$
For any point $x=(x', s)\in \Om$, note that $|\p_nu_t(x',
\rho(x'))|\le CK|x'|$, where $|x'|\le Cs^{1/2}$ by the uniform
convexity of $\pom$. Hence similarly we have $|\p_nu_{t'}(x',
s)|\le CKs^{1/2}$. It follows that
$$|\p_nu_{t'}(x)| \le CKs^{1/2}
        \ \ \ \text{on}\ \ \{x_n=s\}. \tag 3.24$$

Denote $T_i=\p_i+\sum_{j<n} \rho_{x_ix_j}(0)(x_j\p_n-x_n\p_j)$ and
let
$$z(x)=\pm T_i(u_{t'}-\phi) + B (|x'|^2+ s^{-1}x_n^2)-\wtt C Kx_n. $$
By differentiating equation (1.1) with respect to $T_i$, one has
[8]
$$\Cal L z=\pm [T_i (\log f_{t'}) - \Cal L(T_i\phi)]
     +2B (\sum_{i<n-1} u_{t'}^{ii} +  s^{-1}u_{t'}^{nn}), \tag 3.25$$
where $\Cal L =u_{t'}^{ij} \p_i\p_j$ is the linearized operator of
the equation $\log \det D^2 u_{t'}=\log f_{t'}$ , and
$\{u_{t'}^{ij}\}$ is the inverse of the Hessian matrix
$\{D^2u_{t'}\}$.

Let $G=\Om\cap \{x_n<s\}$. First we verify $z\le 0$ on $\p G$. By
subtracting a smooth function  we may assume that $D\phi(0)=0$. By
the boundary condition we have $|T_i(u_{t'}-\phi)|\le C|x|^2$ on
$\pom\cap\p G$. Hence for any given $B>0$, we may choose $\wtt C$
large such that $z\le 0$ on $\p G\cap \pom$. On the part $\p
G\cap\{x_n=s\}$, by (3.23) and (3.24),
$$|T_i(u_{t'}-\phi)|(x)\le CKs+|x'|\,|\p_n u_{t'}|\le CKs. $$
Hence we have $z\le 0$ on $\p G$.

Next we verify that $\Cal L z\ge 0$ in $G$. We compute
$$|D\log f_{t'}|\le C{\tau'}^{\alpha-1}
                  \le C{t'}^{\eps_0 (\alpha-1)},\tag 3.26$$
where $\tau'={t'}^{\eps_0}$ ($\eps_0=1/4n$)  as in (2.3). Observe
that
$$\sum_{i<n} u_{t'}^{ii}+s^{-1}u_{t'}^{nn}\ge
   ns^{-1/n} [\det D^2 u_{t'}]^{-1/n} \ge C s^{-1/n}.$$
Hence we may choose the constant $B$ large, independent of $K, t,
t'$, such that $\Cal L z\ge 0$ in $G$. Now by the maximum
principle we see that $z$ attains its maximum at the origin. It
follows $z_n\le 0$, namely $ |\p_i\p_n u_{t'}(0)|\le CK$.
$\square$

\vskip15pt

Now we choose a fixed small constant $t_0>0$, and for $k=1, 2,
\cdots$, let
$$t_k=t_{k-1}^{1+ \th}=\cdots = t_0^{(1+ \th)^k},
                           \ \  \th=\frac {\alpha}{16n}, \tag 3.27$$
and let $u_k=u_{t_k}$ be the solution of (2.5) with $t=t_k$. Then
we have the estimates
$$\align
\p_\xi^2 u_k &\le C
          \ \ \ \text{in}\ \ \ D_{t_k/8},    \tag 3.28a\\
|\p_\xi\p_\gamma u_k| &\le C^k/\sqrt {t_0}
                \ \ \ \text{in}\ \ D_{t_k/8},  \tag 3.28b\\
\p_\gamma^2 u_k &\le C^k/t_0
           \ \ \ \text{in}\ \ \ D_{t_k/8}.     \tag 3.28c\\
\endalign $$
where the constant $C$ is independent of $k$ and $t_0$. Note that
$$ C^k =O(|\log t_k|^m)\tag 3.29 $$
for some $m>0$ depending only on $C$. Hence for sufficiently large
$k$, (3.13) holds with $\beta<1$ sufficiently close to $1$. Hence
in both Lemmas 3.3 and 3.4, we have
$$|u-u_t|(x)\le Ct^{1+\alpha'/2}\dist(x, \pom) \tag 3.30$$
if $t>0$ is sufficiently small. In particular (3.30) holds for
$u_t=u_{t_k}$ and $u=u_{t_{k+1}}$. From (3.28) and (3.29) we also
have an improvement of (2.20) and (2.21), namely for any small
$\delta>0$,
$$\align
a_h & \le Ch^{(1-\delta)/2},\tag 3.31\\
b_h & \ge Ch^{(1+\delta)/2},\tag 3.32\\
\endalign $$
provided $h$ is sufficiently small, where $C$ is independent of
$h$.

With estimate (3.30), we may introduce the notion of {\it affine
invariant neighborhood} (with respect to the origin). Let $\Ga_i,\
(i=1, 2)$, be two convex hypersurfaces which can be represented as
radial graphs. That is $\Ga_i=\rho_i(x)$ for $x\in S^n$, the unit
sphere (or a subset of $S^n$). We say $\Ga_2$ is in the affine
invariant $\delta$-neighborhood of $\Ga_1$, denoted by
$\Ga_2\subset A_\delta(\Ga_1)$,  if
$$(1-\delta)\rho_2\le \rho_1\le (1+\delta)\rho_2. \tag 3.33$$
If $\Ga_2\subset A_\delta(\Ga_1)$, then $T (\Ga)\subset
A_\delta(T(\pom))$ for any affine transformation $T$ which leaves
the origin invariant, namely $T(x)=T\cdot x$ for some matrix $T$.

Estimate (3.30) gives a control of the shape of the level set
$S_{h, u_k}(0)$ for sufficiently large $k$. Let $h=t_{k+1}^2$, by
convexity and (3.30) we have
$$\align
& |u_k-u|(x)\le Ct_k^{1+\alpha'/2}\dist(x, \pom),\\
& |Du_k|(x)\ge h/|x|\ \ \text{for}\ x\in S_{h, u_k}(0) ,\\
\endalign $$
where we assume that $u_k(0)=0$, $Du_k(0)=0$.  It follows
$$S_{h,u}(0)\subset A_\delta (S_{h,u_k}(0)) \tag 3.34$$
with
$$\align
\delta
 &\le \frac {t_k^{1+\alpha'/2}d_x }{h}
 = t_k^{\alpha'/2-1-2\th} d_x\\
 &\le t_k^{\alpha'/2-1-2\th }t_{k+1}
 = t_k^{\alpha'/2-\th}
 \le t_k^{\alpha'/4} \tag 3.35\\
 \endalign $$
up to a constant $C$.  Note that $|x|$ does not appear in (3.35),
and (3.34) also holds with $u$ replaced by $u_{k+1}$.

\vskip10pt

As a consequence we have an estimate for the shape of the level
set $S_{h, u}(y)$ for any $y\in\pom$. By subtracting a linear
function (which depends on $k$), we assume $u_k(0)=0$ and
$Du_k(0)=0$. By the second inequality of (2.16) we have $S_{h,
u_k}(0)\subset D_{t_k/2}$ for $h= C_0t_k^2$. For simplicity we
assume that $C_0=1$.  We define $a_{h,k}$ and $b_{h,k}$ as in
(2.17) and (2.18) with $u=u_k$.  Let
$$\ol b_{h,k}
   =\sup\{t\ {|}\ \ (0,\cdots, 0, t)\in S_{h, u_k}(0)\}.$$
By Lemma 2.3 and convexity,
$$\ol b_{h,k}\ge \frac {h^{1/2}}{a_{h,k}} b_{h,k}
         \ge  \frac {t_0}{C^{2k}} h^{1/2}. $$
Note that $h^{1/2}=t_k=t_{k-1}^{1+ \th}=\cdots = t_0^{(1+
\th)^k}$.  Consequently for any given $\delta>0$,
$$\ol b_{h,k} \ge C h^{(1+\delta)/2}$$
provided $k$ is sufficiently large, where $C=C(\delta, \th, t_0)$.
Let
$$\ol b_h=\sup\{t\ {|} \ \ (0,\cdots, 0, t)\in S_{h, u}(0)\}.$$
By (3.30) it follows $\ol b_h \ge C h^{(1+\delta)/2}$. Hence
$$u(0, x_n)\le Cx_n^{2/(1+\delta)}\tag 3.36 $$
for $x_n=h^{(1+\delta)/2}$ ($h=t_k^2$). As $k>1$ can be chosen
arbitrary, the above estimate holds for all $x_n>0$ small. By
convexity and the boundary estimates (2.2), we then obtain
$$u(x)\le C|x|^{2/(1+\delta)} \tag 3.37$$
for $x\in\Om$ near the origin. Therefore we have the following
$C^{1, \alpha}$ estimate at the boundary.

\proclaim{Theorem 3.1} Let $u$ be a solution of (1.1) (1.2).
Suppose $\pom, \phi$ and $f$ satisfy the conditions in Theorem
1.1.  Then for any $\hat\alpha\in (0, 1)$, we have the estimate
$$ |u(x)-u(x_0)-Du(x_0)(x-x_0)|\le C|x-x_0|^{1+\hat\alpha} \tag 3.38$$
for any $x\in\Om$ and $x_0\in\pom$, where $C$ depends on
$\hat\alpha$.
\endproclaim

Obviously Theorem 3.1 also holds for $u_t$ for any $t>0$, and the
constant $C$ in (3.38) is independent of $t$. In the next section
we will use a different form of (3.38). That is

\proclaim{Lemma 3.5} Let $u$ satisfy (3.38). Then
$$ |Du(y_0)-Du(y)|\le C|y_0-y|^{\hat\alpha} \tag 3.39 $$
for any $y_0\in\pom$ and $y\in\Om$.
\endproclaim

\noo {\it Proof}. Assume $u(0)=0$, $Du(0)=0$, and $y$ is on the
$x_n$-axis. By convexity we have $\p_\nu u(y) \le \frac 1t
[u(y+t\nu)-u(y)]$ for any unit vector $\nu$ such that $y+t\nu
\in\Om$, where  $t=\frac 12 |y|$. By (3.38), $u(y+t\nu), u(y) \le
Ct^{1+\hat\alpha}$. Hence $\p_\nu u(y)\le Ct^{\hat\alpha}$. It
follows that $|Du(y)-Du(0)|\le C|y|^{\hat\alpha}$. Similarly we
have $|\p_n u(y_0)-\p_n u(0)|\le C|y_0|^\alpha$ for $y_0\in\pom$
near the origin. From the boundary condition, we then infer that
$|D u(y_0)- D u(0)|\le C|y_0|^\alpha$.  Hence (3.39) holds.
$\square$

\vskip30pt

\centerline{\bf \S 4. Continuity estimates for second derivatives}

\vskip10pt

Our passage to $C^2$ estimates at the boundary uses a modulus of
continuity estimate for second derivatives proved by Caffarelli,
Nirenberg, and Spruck in their treatment of the Dirichlet problem
for the Monge-Amp\`ere equation [8, 13].

Let $u_t$ be the solution of (2.5). As before we always suppose
the origin is a boundary point and near the origin $\pom$ is given
by (2.6), and $u_t$ satisfies (2.7).

\proclaim{Lemma 4.1} Suppose $u_t$ satisfies (3.1). Then we have
$$|\p_\xi\p_\gamma u_t(x)-\p_\xi\p_\gamma u_t(0)|
                \le \frac {CK^m}{|\log |x|- \log t|}, \tag 4.1$$
where $m=50$, $x\in \pom$, $|x|\le t/2$.
\endproclaim

\noo{\it Proof}. Although Lemma 4.1 is proved in [8, 13], we
provide an outline here in order to display the polynomial
dependence on the eigenvalue bounds of the coefficients.

Let $v=u_t/t^2$, $y=x/t$. Then $v$ is defined on the set $\{\wtt
\rho (y') < y_n < 1\}$, where $\wtt\rho(y')= \frac 1t\rho(ty')$.
By (2.2), $u_{\xi\xi}\ge C>0$. By the upper bound in (2.16),
$u_t(0, x_n)\ge Cx_n^2$. Hence we have
$$v\ge C\ \ \ \text{on}\ \ \{y_n=1\} \tag 4.2$$
for some positive constant $C$. By (3.1) and Lemma 3.1, we have
$$ C^{-1} K^{-2}\le D^2 v\le CK^2 \ \ \
   \text{in}\ \ G=B_{1/2}(0)\cap \{y_n>\wtt \rho(y')\}, \tag 4.3$$
where the constant $C$ is independent of $K$.

Let $T=\p_i+(\p_i\wtt\rho) \p_n$. Then $T(v-\psi)=T^2(v-\psi)=0$
on $\p G\cap B_{1/2}(0)$, where $\psi(y)=\phi(ty)/t^2$ and $\phi$
is the boundary value in (1.2). By subtracting a smooth function
we may suppose that $D\phi(0)=0$.  Computation as in \S 4 in [8]
shows that
$$\Cal L (T^2 (v-\psi))\ge -CK^8, \tag 4.4$$
where $\Cal L=v^{ij}\p_i\p_j$. Note that the H\"older continuity
of $f_t$ suffices for (4.4), as in the proof of Lemma 2.1.  By
(4.3), the least eigenvalue $\lam$ and the largest eigenvalue
$\Lambda$ of $D^2 v$ satisfy $C^{-1}K^{-2}\le\lam\le\Lambda\le
CK^2$. Hence
$$z=a|y'|^2-by_n^2+cy_n,$$
is an upper barrier of $T^2(v-\psi)$ (in a neighborhood of the
origin) if we choose $a=C_1K^2$, $b= C_2K^{10}$, $c=C_3 K^{10}$
such that $C_3>>C_2>>C_1>0$. It follows that
$$v_{iin}(0)\le CK^{10}.\tag 4.5$$

Let $h=\wtt C K^{10} |y|^2-v_n$. Then
$$|\Cal L h|\le CK^{12}\ \ \ \text{in}\ \ G. \tag 4.6$$
Making the transformation $z'=y'$, $z_n=y_n-\wtt\rho(y')$ to
straighten the boundary $\pom$ near the origin, we may suppose
$G=B_{1/2}^+=B_{1/2}\cap\{y_n>0\}$. By (4.5), $h$ is convex on
$B_{1/2}(0)\cap\{x_n=0\}$ if $\wtt C$ is chosen large. Hence by
the following lemma 4.2, we obtain
$$|\p_\xi\p_\gamma v(y)-\p_\xi\p_\gamma v(0)|
           \le \frac {CK^m}{|\log |y||}\tag 4.9$$
with $m=50$. Scaling back, we obtain  (4.1).

The following Lemma 4.2 is equivalent to Lemma 5.1 in [8].

\proclaim{Lemma 4.2}  Let $h\in C^2(B^+_{1/2})\cap
C^0(B^+_{1/2}\cup T)$ satisfy
$$\Cal Lh=a^{ij} \p_i\p_j h\le \ol f\tag 4.7$$
in $B^+_{1/2}$, where $T=\p B^+_{1/2}\cap\{x_n=0\}$. Let $\lam$
and $\Lambda$ be the least and the largest eigenvalues of the
matrix  $\{a^{ij}\}$. Suppose $h_{{|} T}$ is convex. Then for $x,
y\in T$ near the origin,
$$|\p_i h(x) -\p_i h(y)|
  \le \frac {C}{|\log |x-y||} \frac \Lambda\lam
  \big(\frac {\ol f+\Lambda}{\lam }\big)^3 \sup (|h|+|Dh|),
     \ \ \        i<n.                      \tag 4.8$$
\endproclaim

The main feature of Lemma 4.2, which we used in this paper, is the
polynomial dependence of the modulus of the logarithm continuity
of $\p_i h$ on the eigenvalues of the matrix $\{a_{ij}\}$.
Alternatively we could have used the boundary H\"older estimate of
Krylov [16], which would imply (4.1) with some modulus of
continuity.

\vskip30pt

\centerline{\bf \S 5. Mixed derivative estimates at the boundary
continued}

\vskip10pt

To prove the $C^{2, \alpha}$ estimates at the boundary, we need a
refinement of Lemma 3.4. Let $t_k$ be as in (3.27) and $u_k$ be
the solution of (2.5) with $t=t_k$.

\proclaim {Lemma 5.1} For any  given small $\sigma>0$, there
exists $K>1$ sufficiently large such that if
$$|\p_\xi\p_\gamma u_k|\le K
          \ \ \ \text{on} \ \ \pom, \tag 5.1$$
then
$$|\p_\xi\p_\gamma u_{k+1}|
   \le (1+\sigma)K \ \ \ \text{on}\ \ \pom ,\tag 5.2$$
where $\xi$ is any unit tangential vector on $\pom$, and $\gamma$
is the unit outward normal to $\pom$.
\endproclaim

The constant $\sigma>0$ will be chosen small enough so that
$$(1+ 10 \sigma)^m\le 1+\frac 12 \th , \tag 5.3$$
where $m=50$ as in (4.1) and $\th=\alpha/16n$ as defined before
Lemma 3.4. We also assume $K$ is sufficiently large and $t_k$
sufficiently small such that
$$\align
K\sigma^2 & >1,\tag 5.4\\
K^{20}t_k & \le \sigma^2. \tag 5.5\\
\endalign $$
Note that (5.5) is satisfied when $k$ is large, see (3.29).
Therefore we can also choose $t_0$ sufficiently small such that
(5.5) holds for all $k$.

\vskip10pt

\noo{\it Proof of Lemma 5.1}. \ The proof is also a refinement of
that of Lemma 3.4. As before we suppose the origin is a boundary
point and near the origin $\pom$ is given by (2.6), and $u_k$
satisfies (2.7). Then by (3.30),
$$|Du_{k+1}|(0)=O(t_k^{1+\alpha'/2})=o(t_{k+1}).\tag 5.6$$
By subtracting a smooth function we assume that $\phi(0)=0$,
$D\phi(0)=0$.

Let $\Cal L=u_{k+1}^{ij}\p_i\p_j$ be the linearized operator of
the equation $\log \det D^2 u_{k+1}=\log f_{t_{k+1}}$. Let
$G=D_{t_{k+1}/8}\cap \{x_n<s\}$, where $s=t_{k+1}^{1/4}$. Let
$$\align
 T &=T_i=\p_i+\sum_{j<n} \rho_{x_ix_j}(0)(x_j\p_n-x_n\p_j)\\
 z(x)& =\pm T_i(u_{k+1}-\phi)
               + \frac 12 (|x'|^2+ s^{-1}x_n^2)-(1+8\sigma)Kx_n.\\
\endalign $$
If $\Cal L z\ge 0$ in $G$ and $z\le 0$ on $\p G$, then by the
maximum principle, $z$ attains its maximum at the origin. Hence
$z_n\le 0$ and so $|\p_i\p_n u_{k+1}(0)|\le (1+10\sigma)K$ if
$\sigma K$ is large enough to control $|D^2 \phi|$. Hence Lemma
5.1 holds. In the following we verify that $\Cal L z\ge 0$ in $G$
and $z\le 0$ on $\p G$.

The verification of $\Cal L z\ge 0$ in $G$ is similar to that in
the proof of Lemma 3.4.  We have
$$\Cal L z=\pm [T (\log f_{t_{k+1}}) - \Cal L(T\phi)]
     +(\sum_{i<n-1} u_{k+1}^{ii} +  s^{-1}u_{k+1}^{nn}).\tag 5.7$$
Similar to (3.26),
$$\align
|T (\log f_{t_{k+1}}) - & \Cal L(T\phi)|
               \le Ct_{k+1}^{\eps_0(\alpha-1)},\\
\sum_{i<n} u_{k+1}^{ii}+s^{-1}u_{k+1}^{nn}
       &\ge ns^{-1/n} [\det D^2 u_{k+1}]^{-1/n} \ge C s^{-1/n},\\
\endalign  $$
where $\eps_0=1/4n$. Hence $\Cal L z\ge 0$ as $s=t_{k+1}^{1/4}$ is
very small.

To verify $z\le 0$ on $\p G$, we divide the boundary $\p G$ into
three parts, that is $\p_1 G=\p G\cap \pom$, $\p_2 G=\p G\cap
\{x_n=s\}$, and $\p_3 G=\p G\cap \p\Om_t$ ($t=t_{k+1}/8$).

First we consider the boundary part $\p_1 G$. For any boundary
point $x\in \pom$ near the origin, let $\xi=\xi_T$ be the
projection of the vector $T=\p_i+\rho_{ij}(0)(x_j\p_n - x_n\p_i)$
on the tangent plane of $\pom$ at $x$. We have
$$ |(T-\xi)|(x) \le C|x|^2.\tag 5.8 $$
Hence for $x\in\pom$ near the origin, we have, by (3.39) and
(5.6), and noting that $\p_\xi(u_{k+1}-\phi)=0$,
$$\align
|T(u_{k+1}-\phi)(x)|
       & \le C |x|^2 |\p_\gamma (u_{k+1}-\phi) (x)|\\
    & \le C|x|^2(|x|^{\hat\alpha}
       +|\p_\gamma (u_{k+1}-\phi) (0)|)\\
    &  \le C|x|^2(|x|^{\hat\alpha}
       +t_{k+1}), \tag 5.9\\
       \endalign $$
where $t_{k+1}=s^4$. Hence $z\le 0$ on $\p_1G$.

Next we consider the part $\p_2 G$. For any given point $x=(x',
s)\in\p_2G$, let $\hat x=(x',\rho(x'))\in\pom$. As above let $\xi$
be the projection of $T(\hat x)$ on $\pom$. Then
$$\p_\xi (u_{k+1}-\phi) (x)
   =\p_\xi (u_{k+1}-\phi) (\hat x)
      +\p_n\p_\xi (u_{k+1}-\phi)(x', s')(s-\rho(x')) $$
for some $s'\in (\rho(x'), s)$.  By Lemma 3.4,
$$|\p_n\p_\xi u_{k+1}|\le
      |\p_\gamma\p_\xi u_{k+1}|+|\p_\xi^2u_{k+1}|\le CK. $$
Note that $\p_\xi (u_{k+1}-\phi) (\hat x)=0$ and $|s-\rho(x')|\le
(1+C|x'|^2)t_{k+1}=2s^4$. Hence by (5.8),
$$\align
|T(u_{k+1}-\phi) (x)|
  &\le |\p_\xi(u_{k+1}-\phi) (x)|
     + |T-\xi|\, |\p_\gamma (u_{k+1}-\phi) (x)|\\
    & \le Cs^4K+C|x|^{2+\hat\alpha}.  \tag 5.10\\
    \endalign $$
where we have used that $|T(x)-\xi|\le |T(x)-T(\hat x)|+|T(\hat
x)-\xi|$ and
$$|T(x)-T(\hat x)|=|\sum_j\rho_{ij}(0)(x_n-\hat x_n)\p_j|
              \le Ct_{k+1}=Cs^4.$$
Hence $z\le 0$ on $\p_2 G$.

Finally we consider the part $\p_3G$. We introduce a mapping
$\eta=\eta_k$ from $\p\Om$ to $\p\Om_t$ for $t = t_{k+1}/8$. For
any boundary point $y\in\pom$, by the strict convexity of $u_k$,
the infimum
$$\inf\{u_k(x)-u_k(y)-Du_k(y)(x-y)\ {|}\ \ x\in\pom_t\}$$
is attained at a (unique) point $z\in\p\Om_t$. We define $\eta
(y)=z$. In other words, $z$ is the unique point in $\pom_t\cap
S_{h, u_k}(y)$ with $h>0$ the largest constant such that $S^0_{h,
u_k}(y) \subset D_t$. The mapping $\eta$ is continuous and one to
one by the strict convexity and smoothness of $\pom_t$. The
purpose of introducing the mapping $\eta$ is to give a more
accurate estimate for $|T(u_k-\phi)|(p)$ for $p\in \pom_t$.

First we consider the point $p=(p_1, \cdots, p_n) \in\pom_t$ such
that $\eta^{-1}(p)$ is the origin. Suppose as before that locally
near the origin, $\pom$ is given by (2.6) and $u_k(0)=0$,
$Du_k(0)=0$. Then $h=\inf_{\pom_t} u_k$. By a rotation of the
coordinates $x'$, we suppose that $\{\p_{ij}
u_k(0)\}_{i,j=1}^{n-1}$ is diagonal. We want to prove that
$$\align
|p_i| & \le \frac {1+4\sigma}{\p_i^2u_k(0)}Kt
             \ \ \ \forall\ i=1, \cdots, n-1, \tag 5.11\\
p_n & \le t+o(t) .\tag 5.12\\
\endalign $$

By (2.2), $\p_i^2 u_k(0)$ has positive upper and lower bounds. By
(3.39), the tangential second derivatives of $u_k$ are H\"older
continuous. Indeed, by the boundary condition $u_k=\phi$ on
$\pom$, we have
$$\p^2_{\xi\zeta}u_k+ \rho_{\xi\zeta}\p_\gamma u_k
    =\p^2_{\xi\zeta}\phi+ \rho_{\xi\zeta} \p_\gamma\phi, \tag 5.13$$
where $\xi$ and $\zeta$ are unit tangential vectors, and $\gamma$
is the unit outer normal. By (3.39), $\p_\gamma u_k$ is H\"older
continuous. Hence
$$|\p^2_{\xi\zeta}u_k(x)-\p^2_{\xi\zeta}u(0)|\le \sigma^2 \tag 5.14$$
for any $x\in\pom$ near the origin and any unit tangential vectors
$\xi$ and $\zeta$.

We will prove (5.11) for $i=1$. By restricting to the 2-plane
determined by the $x_1$-axis and $x_n$-axis, without loss of
generality we may assume that $n=2$. Denote
$$\align
a_h & =\sup\{|x_1|\ {|}\ x\in S_{h, u_k}(0)\}, \\
b_h & =\sup\{x_n\ {|}\ x\in S_{h, u_k}(0)\}. \\
\endalign $$
where $h=\inf_{\pom_t} u_k$. Then it suffices to prove
$$\align
a_h & \le \frac {1+4\sigma}{\p_1^2u_k(0)}Kt ,\tag 5.11$'$\\
b_h & \le t+o(t) .\tag 5.12$'$\\
\endalign $$
Note that we have now $x=(x_1, x_n)$, and the domains $D_t, \Om_t$
denote the restriction on the 2-plane.

Assume the supremum $a_h$ is achieved at $x_h=(a_h, c_h)$. In the
two dimensional case, the level set $\ell:=S_{h, u_k}$ is a curve
in $\bom$, which has an endpoint $\hat x=(\hat x_1, \hat x_n) \in
\pom$ with $\hat x_1>0$.

If $a_h\le Ch^{1/2}$ for some $C>0$ under control, by (2.16) we
have $b_h\ge C_1h^{1/2}$. In this case we have $t\ge C_2h^{1/2}$.
Hence (5.11$'$) holds for sufficiently large $K$.

If $a_h\ge Ch^{1/2}$ (let us choose $C=\sigma^{-2}$), let $\xi,
\zeta, \th_1, \th_2$ be as in the proof of Lemma 2.3. Then
$\th_1+\th_2<\pi/2$. By (5.1) and (5.14),
$$\align
|\p_\gamma u_k (\hat x)| & \le (1+\sigma) K |\hat x|, \tag 5.15\\
|\p_\xi u_k (\hat x)| &
             \ge (1-\sigma) \p^2_1 u_k(0)\, |\hat x|.\tag 5.16\\
\endalign $$
Hence $\text{tg}\th_1\ge \frac
{(1-\sigma)\p_1^2u_k(0)}{(1+\sigma)K}$. Note that
$\text{tg}(\th_1+\th_2)\le c_h/(a_h-\hat x_1)$ by the convexity of
$\ell$. We obtain
$$a_h \le \hat x_1 + \frac {1+2\sigma}{\p_1^2u_k(0)} Kc_h. $$
Recall that $h^{1/2} \le \sigma^2 a_h$ by assumption, and $\hat
x_1\le Ch^{1/2}$ by (2.2).  Hence we obtain
$$a_h \le \frac {1+3\sigma}{\p_1^2u_k(0)} Kc_h. \tag 5.17$$

Suppose $\pom_t$ is locally given by
$$x_n=\rho_t(x').\tag 5.18$$
Then $\rho_t$ is smooth and uniformly convex. It is easy to see
that $\rho_t(0)=t$ and $|D\rho_t|(0)=o(t)$. Hence we have
$$c_h\le t+C_1 a_h^2+o(t) a_h. \tag 5.19$$
By (3.31), $a_h\le Ch^{(1-\delta)/2}$. By (3.36), $h\le
Ct^{2/(1+\delta)}$, where $\delta>0$ can be arbitrarily small as
long as $t$ is sufficiently small. Hence we have $c_h\le t+o(t)$.
Therefore (5.11) holds.

To prove (5.12), assume that the supremum $b_h$ is attained at
$\hat x_h=(d_h, b_h)$. Then $b_h\le \rho_t(d_h)$. Hence
$$b_h\le t+C_1 d_h^2+o(t) d_h\le t+o(t). \tag 5.20$$
Recall that $d_h\le a_h\le Ch^{(1-\delta)/2}$, and by our
definition of $h$, $b_h\ge t$. Hence (5.12) holds.

Now we prove
$$|T(u_k-\phi)|(p)\le (1+6\sigma)K p_n  \tag 5.21$$
at $p=\eta(0)$. Let $\xi$ be the projection of $T(p)$ on the
tangent plane of $\pom_t$ at $p$. We have
$$\align
 |T(p)| & \le 1+ C(p_n+|p|^2) ,\tag 5.22\\
 |(T-\xi)(p)| & \le C(p_n+|p|^2) . \tag 5.23\\
 \endalign $$
Hence
$$|T(u_k-\phi)(p)|
 \le |\p_\xi (u_k-\phi)(p)|
            + C(p_n+|p|^2) |D(u_k-\phi)(p)|.\tag 5.24$$
By (3.39),
$$|D(u_k-\phi)(p)|\le C|p|^{\hat\alpha}. $$
Hence the second term in (5.24) is small. By (5.13), we have
$\p^2_{ij}\phi(0)=\p^2_{ij}u_k(0)$ for $i, j=1, \cdots, n-1$
(recall that we assume $D\phi(0)=0$ at the beginning). Hence near
the origin we have, by the Taylor expansion and (5.11),
$$\align
|\p_i\phi(p)| &\le (1+\sigma)|p_j\,\p_i\p_ju_k(0)|\\
   &\le (1+5\sigma)Kp_n .\tag 5.25\\
   \endalign$$
By our definition of the mapping $\eta$, $\p_\xi u_k=0$ at $p$.
(This is the purpose of introducing the mapping $\eta$). Hence
$$|\p_\xi (u_k-\phi)(p)|\le (1+6\sigma)Kp_n.\tag 5.26$$
By (5.24) we therefore obtain (5.21).

Next we prove (5.21) for any given $p\in\p_3G$. Let
$y=\eta^{-1}(p)$, where $\eta$ is the mapping introduced above.
Then by (5.14) we have, similarly to (5.11),
$$|p_i-y_i|\le \frac {1+5\sigma}{\p_i^2 u_k(0)}Kt. \tag 5.27$$
Choose a new coordinate system such that $y$ is the origin and the
positive $x_n$-axis is the inner normal at $y$. Subtract a linear
function from both $u_k$ and $\phi$ (which does not change the
value of $T(u_k-\phi)$) such that $Du_k(y)=0$. As above let $\xi$
be the projection of $T(p)$ on the tangent plane of $\pom_t$ at
$p$. By (3.39), $|Du_k|, |D\phi|\le \sigma^2$ in $G$. Hence
$$\align
 |Tu_k(p)| & \le |\p_\xi u_k(p)|+|T(p)-\xi|\, |Du_k(p)|\le Cp_n,\\
|T\phi(p)| & \le |\p_\xi \phi (p)|+ Cp_n .\\
\endalign $$
By (5.13) and noting that $|D\phi|\le \sigma^2$, we have, similar
to (5.14),
$$|\p^2_{\xi\zeta}\phi(x)-\p^2_{\xi\zeta}u_k(0)|\le \sigma^2.$$
Hence as (5.25) we have
$$|\p_\xi\phi(p)|\le (1+6\sigma)Kp_n. $$
Hence (5.21) holds at any point $p\in\p_3G$.

With (5.21) we are now in position to prove $z\le 0$ on $\p_3G$.
By (3.30),
$$|u_{k+1}-u_k|(x)\le Ct_k^{1+\alpha'/2} t
                  \ \ \ x\in\pom_t. $$
Hence by (3.28a),
$$|\p_\xi (u_{k+1}-u_k)(x)|\le C (t_k^{1+\alpha'/2} t)^{1/2}
     \le C t_k^{\alpha'/8} t \ \ \ x\in\pom_t,$$
where $\xi$ is any unit tangential vector to $\pom_t$. Hence
$$\align |T (u_{k+1}-u_k)(x)|
  & \le |\p_\xi (u_{k+1}-u_k)|+ C(t+|x|^2) |D(u_k-\phi)| \\
  & \le C t_k^{\alpha'/8} t+ Cx_n .\\
  \endalign $$
In view of (5.21), it follows that
$$|T(u_{k+1}-\phi)(x)|
       \le (1+7\sigma)K x_n\ \ \ x\in\pom_t.\tag 5.28$$
From (5.28) and noting that $\sigma K>>1$, we obtain $z\le 0$ on
$\p_3 G$. This completes the proof. $\square$

By Lemma 5.1, we improve (3.28) to
$$\align
\p_\xi^2 u_k & \le C
           \ \ \ \text{in}\ \ D_{t_k/8}  ,  \tag 5.29a\\
|\p_\xi\p_\gamma u_k| &\le C(1+\sigma)^k
           \ \ \ \text{in}\ \ D_{t_k/8}  ,   \tag 5.29b\\
\p_\gamma^2 u_k &\le C (1+\sigma)^{2k}
           \ \ \ \text{in}\ \ D_{t_k/8}  ,    \tag 5.29c\\
\endalign $$
where $C$ depends only on $n, \pom, f, t_0$, and $\phi$.

Now we apply the estimate (4.1) to the section $S^0_{h, u_k}(0)$,
where
$$h=t_{k+1}^2=t_k^{2(1+ \th)}, \ \ \  \th=\alpha/16n.$$
For any $x \in \pom \cap S^0_{h, u_k}$, we have by (2.2),
$$|x| \le Ch^{1/2} \le Ct_{k+1} . $$
By (4.1),
$$|\p_\xi\p_\gamma u_k(x) - \p_\xi\p_\gamma u_k(0)| \le
 \frac {[C(1+\sigma)^k ]^m }{|\log |x|-\log t_k|}$$
By our definition, $t_k=t_{k-1}^{1+ \th}=\cdots =t_0^{(1+
\th)^k}$. We obtain, by the choice of $\sigma$ in (5.3),
$$|\p_\xi\p_\gamma u_k(x) - \p_\xi\p_\gamma u_k(0)|
 \le \wtt C\frac {(1+ \th/2)^k}{(1+ \th)^k} ,\tag 5.30$$
where $\wtt C$ depends only on $n, \pom, f, \phi$ and $t_0$, and
is independent of $k$.

\vskip10pt

\noo{\it Proof of Theorem 1.1}. We will first prove
$$\sup_{x\in\Om} |D^2 u(x)|\le C. \tag 5.31$$
Suppose the origin is a boundary point such that
$\Om\subset\{x_n>0\}$. We will prove $D^2 u$ is bounded at the
origin. By making a linear transformation of the form
$$\align
y_n & =x_n\\
y_i & =x_i-\alpha_i x_n,\ \ \ i=1, \cdots, n-1, \tag 5.32\\
\endalign $$
we may suppose $\p_i\p_n u_k(0)=0$, where by (5.29b),
$$ |\alpha_i|\le C(1+\sigma)^k\le C|\log h|.$$
Hence the boundary part $\{x\in\pom\ {|}\ u_k(x)<h\}$ is smooth
and uniformly convex after the transformation (5.32).  By (5.30)
there is a sufficiently large $k_0$ such that when $k\ge k_0$,
$$|\p_\xi\p_\gamma u_k(x)|\le C \tag 5.33$$
for $x\in\pom$ with $|x|<t_{k+1}$. Hence from (2.20) and (2.21),
$$\align
a_{h, k} & =\sup\{|x'|\ {|}\ x\in S_{h, u_k}(0)\}
                            \le \wtt C h^{1/2}, \\
b_{h, k} & =\sup\{x_n\ {|}\ x\in S_{h, u_k}(0)\}
                 \ge h^{1/2}/\wtt C \tag 5.34  \\
\endalign $$
for some $\wtt C>0$ depending only on $n, f, \phi$ and $\pom$, but
independent of $k$. That is the section $S^0_{h, u_k}$ has a good
shape, as defined in (2.24).

By (3.34), $S^0_{h, u}$ also has a good shape for $h\le
t_{k+1}^2$. Now the perturbation argument [4], see Section 6,
implies that
$$C_1|x|^2\le u(x)\le C_2|x|^2, \tag 5.35$$
where we assume $u(0)=0$, $Du(0)=0$. Furthermore, $|D^2 u(x)| \le
C$, for $x\in \Om$ near the origin. Making the inverse
transformation of (5.32), we obtain (5.31) for $x$ near the
origin. The interior second order derivative estimate was
established in [4]. Hence (5.31) holds.

Estimate (5.31) implies the Monge-Amp\`ere equation is uniformly
elliptic, and hence the $C^{2, \alpha}$ estimate follows [2,19].
$\square$

\vskip8pt

\noo{\bf Remark}. Estimate (5.30) actually implies a continuity
estimate for the mixed second derivatives of $u$ on the boundary.
By the $C^{1, \alpha}$ estimate (Lemma 3.5) and the equation
itself, we can then infer a continuity estimate for $D^2 u$ on the
boundary. However, unless the inhomogeneous term $f$ is smoother,
we shall need to use the perturbation argument of the next section
to derive continuity estimates for $D^2 u$ near the boundary.

\vskip30pt

\centerline{\bf \S 6. The perturbation argument}

\vskip10pt

In this section we provide the perturbation argument [4] which
enables us to proceed from a level set of good shape to second
derivative estimates.

\proclaim{Theorem 6.1} Let $u$ be a convex solution to (1.1)
(1.2). Suppose there is an $h_0>0$ such that for any boundary
point $y\in\pom$, $S_{h_0, u}^0(y)$ has a good shape. Then under
the assumptions of Theorem 1.1, $u$ is $C^{2, \alpha}$ smooth up
to the boundary.
\endproclaim

\noo{\it Proof}. Let the origin be a boundary point such that
$\Om\subset \{x_n>0\}$. By subtracting a linear function we
suppose
$$u(0)=0,\ \ \ Du(0)=0. \tag 6.1$$
By a rescaling $u\to u/h_0$, $x\to x/\sqrt{h_0}$, we may suppose
$h_0=1$ and
$$|f(x)-f(0)|\le \eps |x|^\alpha   \tag 6.2$$
for some $\eps>0$ sufficiently small. For simplicity we suppose
$f(0)=1$. By (2.2) we have
$$C^{-1}\le u_{\xi\xi} \le C\ \ \ \text{on}\ \ \pom\tag 6.3$$
for any unit tangential vector $\xi$. First we need two lemmas.

\proclaim{Lemma 6.1}\ \ Let $u_i$, $i=1,2$, be two convex
solutions of $\det D^2 u=1$ such that $u_1=u_2$ on $\pom$. Suppose
$\|u_i\|_{C^{2, \alpha}}\le C_0$ in $S^0_{1, u_1}(0)$. Then if
$$|u_1-u_2|\le \delta\ \ \ \text{in}\ \ S^0_{1, u_1}\tag 6.4$$
for some sufficiently small $\delta>0$, we have
$$ |D^2 (u_1-u_2)|\le C\delta
              \ \ \ \text{in}\ \ S^0_{1/2, u_1}. \tag 6.5$$
\endproclaim

\noo{\it Proof.}\ \ We have
$$\align
\det D^2 u_2-\det D^2 u_1
    & = \int_0^1 \frac{d}{dt}\det[D^2 u_1+t(D^2 u_2-D^2 u_1)]dt\\
    & =a_{ij}(x) \p_i\p_j(u_2-u_1) = 0,\tag 6.6\\
\endalign $$
where $L=a_{ij}(x)\p_i\p_j$ is a linear, uniformly elliptic
operator with H\"older continuous coefficients. By the Schauder
estimates for linear elliptic equations, we obtain (6.5).
$\square$

\proclaim{Lemma 6.2}\ \ Let $u$ be as above such that $S^0_{1, u}$
has a good shape. Then for $h\in (0, 1/4]$,
$$S_{h, u}\subset N_\delta (h^{1/2} E)\tag 6.7$$
with
$$\delta \le C(h^{(1+\hat\alpha)/2}+h^{-1/2}\eps) ,\tag 6.8$$
where $\hat\alpha$ is any constant in $(0, 1)$, $N_\delta$ denotes
the $\delta$-neighborhood, $E$ is an ellipsoid of good shape.
\endproclaim

\noo{\it Proof.}\ \ Let $v$ be the solution of
$$\det D^2 v=f(0)=1\ \ \ \text{in}\ \ S^0_{1,u}$$
such that $v=u$ on $\p S^0_{1, u}$. Since $u=\phi\in C^3$ on
$\pom$ and $\pom\in C^3$, from [22] we have $v\in C^{2,
\hat\alpha}(S^0_{3/4, u})$ $\forall\ \hat\alpha\in (0, 1)$. By the
Taylor expansion,
$$v(x)=v(0)+v_i(0)x_i+\frac 12 v_{ij}(0)x_ix_j
                                +O(|x|^{2+\hat\alpha}),$$
we have, on $S_{h, v}(0)$,
$$ C^{-1} h^{1/2}\le |Dv|\le C h^{1/2}.\tag 6.9$$
Hence
$$S_{h,v}(0)\le N_{\hat \delta} (h^{1/2} E)$$
with $\hat \delta\le C h^{(1+\hat\alpha)/2}$, where $E$ is the
ellipsoid $\{x\in\R^n\ {|}\ \frac 12 v_{ij}(0)x_ix_j=1\}$.

By (6.2) it is easy to verify that $|u-v|\le C\eps$, and by (6.3)
we have $|Dv(0)|\le C\eps$. Hence by (6.9), we have
$$S^0_{h-Ch^{-1/2}\eps, v}(0)\le S^0_{h,u}(0)
     \le S^0_{h+Ch^{-1/2}\eps, v}(0) \tag 6.10$$
provided $\eps<<h^{1/2}$. Hence
$$S_{h, u}\subset N_{Ch^{-1/2}\eps}(S_{h, v})
    \subset N_{Ch^{(1+\hat\alpha)/2}+Ch^{-1/2}\eps}(h^{1/2}E).
                                     \eqno \square $$

\noo{\it Proof of Theorem 6.1 continued}:\ \ \ Let $u_k$, $k=0, 1,
\cdots$, be the solution of
$$\align
\det D^2 u_k & = 1\ \ \ \text{in}\ \ S^0_{4^{-k}, u},\\
         u_k & = u\ \ \ \text{on}\ \ \p S^0_{4^{-k}, u}.\\
\endalign $$
Since $S^0_{1, u}$ has a good shape, by the regularity of the
Monge-Amp\`ere equation, we have $\|u_0\|_{C^{2, \alpha}
(S^0_{3/4,u})}\le C $. Denote
$$\omega_k=\sup\{|f(x)-1|\ {|}\ x\in S^0_{4^{-k}, u}\},$$
where $f(0)=1$ by assumption. By the comparison principle we have
$|u-u_0|\le C\omega_0$.  Hence if the constant $\eps$ in (6.2) is
sufficiently small, $S^0_{1/4, u}$ has a good shape. It follows
$\|u_1\|_{C^{2, \hat\alpha} (S^0_{3/16, u_1})}\le C$. Note that
$|u_1-u_0|\le C\omega_0$. By Lemma 6.1 we obtain
$$|D^2 u_0(x)-D^2 u_1(x)|\le C\omega_0
     \ \ \text{for}\ \ x\in S^0_{4^{-2}, u_1}. \tag 6.11$$
It follows that $2^2 S^0_{4^{-2}, u_1}$ has a good shape, where
$t\Om=\{x\in\R^n\ {|}\ \ t x\in \Om\}$.

Let $R_k=\sup\{|x|\ {|}\ \ x\in S^0_{4^{-k}, u}\}$, namely
$B_{R_k}(0)$ is the smallest ball containing $S^0_{4^{-k}, u}$.
By (6.11) there is a constant $\beta>0$ such that
$$R_1<(1-\beta)R_0. \tag 6.12$$

For $k=1, 2, \cdots$, applying the same argument to $\hat u_0:=
4^ku_k(2^{-k}x)$ and $\hat u_1:= 4^ku_{k+1}(2^{-k}x)$, we obtain
$$|D^2 u_k(x)-D^2 u_{k+1}(x)|\le C\omega_k
     \ \ \text{for}\ \ x\in S^0_{4^{-k-2}, u_{k+1}}. \tag 6.13$$

From (6.2) and by induction we have
$$\align
R_k & \le (1-\beta)R_{k-1}\le C(1-\beta)^k ,\\
\omega_k & \le C\eps (1-\beta)^{\alpha k}.\\
\endalign $$
Hence we obtain from (6.13),
$$|D^2 u_0(x)-D^2 u_{k+1}(x)|\le C\sum_{i=0}^k\omega_i
     \ \ \text{for}\ \ x\in S^0_{4^{-k-2}, u_{k+1}}, \tag 6.14$$
where the right hand side $\le C\eps$. Hence $S^0_{4^{-k},
u}=S^0_{4^{-k}, u_k}$ has a good shape. From (6.14) we see that
$\{D^2 u_{k+1}(0)\}$ is convergent. Hence $u$ is twice
differentiable at $0$, and $D^2 u(0)=\lim_{k\to\infty} D^2
u_k(0)$. Moreover, $(D^2 u)$ is positive definite, so the
Monge-Amp\`ere equation (1.1) is uniformly elliptic. The H\"older
continuity of $D^2 u$ follows from [2, 19].

The H\"older continuity of $D^2 u$ also follows from (6.14)
immediately. Indeed, let $\hat x$ be a point in $\bom$ near the
origin. Choose $k_0$ such that $\hat x\in S_{4^{-k_0-1}, u} (0)$.
For $k\ge k_0$, let $\hat u_k$ be the solution of
$$\align
\det D^2 \hat u_k & = \hat f_k
                \ \ \ \text{in}\ \ S^0_{4^{-k}, u}(\hat x),\\
 \hat u_k & = u\ \ \ \text{on}\ \ \p S^0_{4^{-k}, u}(\hat x),\\
\endalign $$
where $\hat f_k=\inf\{f(x)\ {|}\ x\in S^0_{4^{-k}, u}(\hat x)\}$.
Then similarly we have
$$|D^2 \hat u_{k_0}(\hat x)-D^2 \hat u_{k+1}(\hat x)|
   \le C\sum_{i=k_0}^k\hat \omega_i,    \tag 6.15$$
where $\hat\omega_k\le \sup\{|f(x)-f(\hat x)|
   \ {|}\ x\in S^0_{4^{-k}, u}(\hat x)\}$.
Since $f$ is H\"older continuous, $\sum_{i=k_0}^\infty\hat
\omega_i\le Cd_0^\alpha$ and $\sum_{i=k_0}^\infty\omega_i\le C
d_0^\alpha$, where $d_0$ is the diameter of the set
$S_{4^{-k_0-1}, u} (0)$. From (6.14), (6.15), and the interior
smoothness of $u_{k_0}$, and by choosing appropriate $k_0$, we
obtain the H\"older continuity at the origin,
$$|D^2 u(\hat x)-D^2 u(0)|
  \le C|\hat x|^{\alpha'}\tag 6.16$$
for some $\alpha'\in (0, \alpha)$. From (6.16) we obtain the
global H\"older continuity for $D^2 u$. Indeed, let $x, y\in\Om$
and close to $\pom$. If $|x-y|\ge \delta_0 (\dist(x,
\pom)+\dist(y, \pom))$ for some constant $\delta_0>0$, let $\hat
x, \hat y\in\pom$ be the boundary points closest to $x, y$. Then
by (6.16) (denote $\Cal A(x, y)=|D^2 u(x)-D^2 u(y)|$ for short)
$$\Cal A(x, y)\le \Cal A(x, \hat x)+\Cal A(\hat x, \hat y)
  +\Cal A(\hat y, y)\le C|x-y|^{\alpha'}.$$
Otherwise the estimate for $\Cal A(x, y)$ is equivalent to the
interior one [4]. $\square$

\noo{\it Remark 6.1}. For the estimate (6.16), if $\hat x$ is also
a boundary point, the proof uses only the H\"older continuity of
$f$ in the sets $S^0_{h, u}(x)$ for $x\in\pom$. Hence if $f$
satisfies (3.16), $D^2 u$ is H\"older continuous on $\pom$. We do
not require that $f$ is H\"older in $\Om$.

\noo{\it Remark 6.2}.  We have actually proved that $D^2 u$ is
continuous if $f$ is Dini continuous, that is if
$$\int_0^1\frac {\omega(t)}{t} dt<\infty, $$
where $\omega(t)=\sup\{|f(x)-f(y)|\ {|}\ \ |x-y|<t\}$, so that the
right hand side of (6.14) is convergent.

\noo{\it Remark 6.3}.  For the interior $C^{2, \alpha}$ estimate,
the condition that $S^0_{h_0, u}$ has a good shape is
automatically satisfied if $u$ is a strictly convex solution,
since the convex set $S^0_{h_0, u}$ can be normalized by a linear
transformation. However for the $C^{2, \alpha}$ estimate at the
boundary, we can only do linear transformation of the form (5.32)
with relatively small $\alpha_i$, and must prove (5.34) for $u$ so
that the level set has a good shape. Other linear transformations
may worsen the boundary condition.

\vskip30pt

\centerline{\bf \S 7. Application to the affine mean curvature
equation}

\vskip10pt

In this section we prove Theorem 1.2.  First we prove the
uniqueness of solutions.

\proclaim{Lemma 7.1} There is at most one uniformly convex
solution $u\in C^4(\Om)\cap C^2(\bom)$ of the second boundary
value problem (1.4)-(1.6).
\endproclaim

\noo{\it Proof}.\ \ Suppose both $u_1$ and $u_2$ are solutions. We
have, by the concavity of the affine area functional $A$,
$$\align
A(u_1)-A(u_2)
 &= \int_\Om \bigg (\det D^2 u_1)^{1/(n+2)}
              -(\det D^2 u_2)^{1/(n+2)}\bigg)\\
 &\le \frac {1}{n+2}\int_\Om w_2U_2^{ij} D_{ij}(u_1-u_2)\\
 &= \frac {1}{n+2}\big[\int_\pom \gamma_i D_j(u_1-u_2) w_2 U_2^{ij}
 + \int_\Om (u_1-u_2) f(x) \big].    \\
\endalign $$
where we have used the divergence free relation $\sum_i \p_i
U^{ij}=0$ $\forall\ j$. Similarly we have
$$A(u_2)-A(u_1) \le
  \frac {1}{n+2}\big[\int_\pom \gamma_i D_j(u_2-u_1) w_1 U_1^{ij}
   - \int_\Om (u_1-u_2) f(x) \big].$$
Note that $w_1=w_2$ on $\pom$. Hence
$$0 \le \int_\pom  w_1\gamma_i D_j(u_1-u_2) (U_2^{ij}-U_1^{ij})
 = -\int_\pom  w_1\gamma_i D_j(u_1-u_2) (U_1^{ij}-U_2^{ij}). $$
For any given boundary point, suppose $e_n=(0, \cdots, 0, 1)$ is
the inner normal there. Then $\gamma=-e_n$, and the right hand
side of the above inequality is equal to
$$-\int_\pom w_1 D_n (u_1-u_2) (U_1^{nn}-U_2^{nn}) ,$$
where $U^{nn}=\det (u_{x_ix_j}){|}_{i, j=1}^{n-1}$. Since
$u_1=u_2$ on $\pom$, we have
$$U_1^{nn}-U_2^{nn}>0\ \ \
  \text{if}\ \ \frac {\p u_1}{\p x_n} < \frac {\p u_2}{\p x_n}.$$
Hence we obtain
$$ 0\le \int_\pom w_1 D_n (u_1-u_2) (U_1^{nn}-U_2^{nn}) <0, $$
which implies $Du_1=Du_2$ on $\pom$. Hence $u_1=u_2$ by the
concavity of the affine area functional. This completes the proof.
$\square$

In the following we always assume that $u\in C^4(\bom)$ is a
uniformly convex solution of (1.4)-(1.6) and the conditions of
Theorem 1.2 hold. By Aleksandrov's maximum principle [13], $u\in
W^{4, p}(\Om)$ ($p\ge n$) suffices for the estimates below. Note
that $u\in W^{4,1}_{loc} (\Om)\cap C^2(\bom)$ suffices for Lemma
7.1. The following lemma is taken from [21]

\proclaim{Lemma 7.2} There exists a constant $C>0$ such that any
solution $u$ of (1.4) satisfies
$$\align
 C^{-1} \le w & \le C\ \ \ \text{in}\ \ \Om, \tag 7.1\\
 |w(x)-w(x_0)|& \le C|x-x_0|
        \ \ \ \forall\ \ x\in\Om, x_0\in\pom, \tag 7.2\\
\endalign$$
where $C$ depends only on $n$, $\diam(\Om)$, $\sup_\Om |f|$, and
$\sup_\Om |u|$.
\endproclaim

\noo{\it Proof.}  Let $z=\log w -u$. If $z$ attains its minimum at
a boundary point, by the boundary condition (1.6) we have $w\ge C$
in $\Om$. Let us suppose $z$ attains its minimum at an interior
point $x_0\in\Om$. At this point we have
$$\align
0 & = z_i=\frac {w_i}{w}-u_i,\\
0 &\le z_{ij} = \frac {w_{ij}}{w}-\frac {w_iw_j}{w^2}-u_{ij}\\
\endalign$$
as a matrix. Hence
$$0\le u^{ij} z_{ij} \le \frac {f}{d^\th} -n$$
where $d=\det D^2 u$, $\th=1/(n+2)$. We obtain $d(x_0)\le C$.
Since $z(x)\ge z(x_0)$, we obtain
$$w(x)\ge w(x_0) \text{exp} (u(x)-u(x_0)). \tag 7.3$$
The first inequality in (7.1) follows.

Next let $z=\log w+A|x|^2$. If $z$ attains its maximum at a
boundary point, by (1.6) we have $w\le C$ and so (7.1) holds. If
$z$ attains its maximum at an interior point $x_0$, we have, at
$x_0$,
$$\align
0 & = z_i=\frac {w_i}{w}+2Ax_i,\\
0 &\ge z_{ii} = \frac {w_{ii}}{w}-\frac {w_i^2}{w^2}+2A .\\
\endalign$$
Suppose $(D^2 u)$ is diagonal at $x_0$. Then
$$0  \ge u^{ij} z_{ij}
  = \frac {f}{d^\th} -4A^2 x_i^2 u^{ii}+2Au^{ii}
 \ge \frac {f}{d^\th} +Au^{ii} \tag 7.4$$
if $A$ is small. Observe that
$$d^\th \sum u^{ii}\ge C(\sum u^{ii})^{2/(n+2)}$$
We obtain $\sum u^{ii}\le C$, and hence (7.1) is proved.

Let $v$ be a smooth, uniformly convex function in $\Om$ such that
$v=\psi$ on $\pom$ and $D^2 v\ge K$. Then
$$U^{ij}v_{ij}
     \ge K \sum U^{ii}\ge CK[\det D^2 v]^{(n-1)/n}\ge CK.$$
Hence if $K$ is large enough, $v$ is a lower barrier of $w$
(regarding (1.4) as a second order elliptic equation of $w$). We
thus obtain
$$w(x)-w(x_0)\ge -C|x-x_0|
        \ \ \ \forall\ \ x\in\Om, x_0\in\pom.\tag 7.5$$
Similarly one can construct an upper barrier for $w$. Hence (7.2)
holds. $\square$

In (7.3) the lower bound for $w$ depends on the uniform estimate
for $u$. To obtain the uniform estimate for $u$, we in turn need
the lower bound for $w$, namely the upper bound for $\det D^2 u$.
To avoid the mutual dependence we assume $f\le 0$, so that $w$
attains its minimum on the boundary by the maximum principle. This
condition can be relaxed to $f\le \eps$ for some $\eps>0$ small
but cannot be removed completely, as is easily seen by solving
equation (1.4) in the one dimensional case.

\proclaim{Lemma 7.3}\ \ Let $u\in C^4(\bom)$ be a solution of the
boundary value problem (1.4)-(1.6). Then we have the estimate
$$\sup_\Om |D^2 u|\le C ,\tag 7.6$$
where $C$ depends only on $n, \pom$, $\|f\|_{L^\infty}$,
$\|\phi\|_{C^4(\bom)}$, $\|\psi\|_{C^4(\bom)}$, and $\inf \psi$.
\endproclaim

\noo{\it Proof.}\ \ Consider the Monge-Amp\`ere equation
$$\det D^2 u=w^{-(n+2)/(n+1)}\ \ \ \text{in}\ \ \Om.\tag 7.7$$
By Lemma 7.2, the right hand side of (7.7) is positive and
satisfies condition (3.16). Hence by the argument in the preceding
sections,  $D^2 u$ is bounded and H\"older continuous on the
boundary, see Remark 6.1. For any $\delta>0$, by (7.1) the
solution of the linearized Monge-Amp\`ere equation
$$U^{ij}w_{ij}  =f \ \ \ \text{in}\ \ \Om\tag 7.8$$
is H\"older continuous [7], namely $\det D^2 u\in
C^\alpha(\Om_\delta)$ for some $\alpha\in (0, 1)$. Hence $u\in
C^{2, \alpha}(\Om_\delta)$ [4]. So we are left to consider a point
$\hat x\in\Om$ near the boundary. Choosing an appropriate
coordinate system, we assume that $\hat x$ is on the positive
$x_n$-axis, the origin is a boundary point, and
$\Om\subset\{x_n>0\}$. Suppose $u(0)=0$, $Du(0)=0$. Then the
argument of the preceding sections apply, with $\th=\frac 1{16n}$,
and we conclude as before the quadratic growth estimate (5.35).
Let $\hat h$ is the largest constant such that $S^0_{\hat h,
u}(\hat x)\subset\Om$. By (5.35), the section $S^0_{\hat h,
u}(\hat x)$ has a good shape. Hence the argument in [7] applies,
and we also conclude that $w$ is bounded and H\"older continuous
near $\hat x$. Hence (7.6) holds. $\square$

\proclaim{Lemma 7.4}\ \ If $f\in L^\infty(\Om)$, then for any
$p>1$, we have
$$\|u\|_{W^{4, p}(\Om)}\le C, \tag 7.9$$
where $C$ depends only on $n, p, \pom$, $\|f\|_{L^\infty}$,
$\|\phi\|_{C^4(\bom)}$, $\|\psi\|_{C^4(\bom)}$,  and $\inf \psi$.
If $f\in C^\alpha(\bom)$, $\phi\in C^{4, \alpha}(\bom)$,  $\psi\in
C^{4, \alpha}(\bom)$, and $\pom\in C^{4, \alpha}$  for some
$\alpha\in (0, 1)$, then
$$u\in C^{4, \alpha}(\bom)\le C\tag 7.10$$
where $C$ depends in addition on $\alpha$.
\endproclaim

\noo{\it Proof}. Regard the fourth order equation (1.4) as a
system of two second order partial differential equations (7.7)
(7.8). By estimate (7.6), both (7.7) and (7.8) are uniformly
elliptic. It follows that $w$ is H\"older continuous up to the
boundary and so $u\in C^{2, \alpha}(\bom)$ [2, 19]. Hence (7.8) is
a linear, uniformly elliptic equation with H\"older coefficients.
Hence $w\in W^{2, p}(\Om)$ for any $p<\infty$. From (7.7) we also
conclude the global $C^{4, \alpha}$ a priori estimate for $u$.
$\square$

\vskip10pt

\noo{\it Proof of Theorem 1.2}.  We have proved the uniqueness and
established the a priori estimate for solutions of (1.4)-(1.6). To
prove the existence of solutions we use the degree theory as
follows.

For any positive $w\in C^{0, 1}(\bom)$, let $u=u_w\in C^{2,
\alpha}(\bom)$ be the solution of
$$\align
\det D^2 u & = w^{-(n+2)/(n+1)}\ \ \text{in}\ \ \Om, \tag 7.11\\
  u & =\phi\ \ \ \text{on}\ \ \pom .\\
  \endalign $$
Next let $w_t$, $t\in [0, 1]$, be the solution of
$$\align
 U^{ij}w_{ij} & = tf(x)\ \ \text{in}\ \ \Om ,\tag 7.12\\
          w_t & =t\psi+(1-t)\ \ \ \text{on}\ \ \pom .\\
\endalign $$
We have thus defined a compact mapping $T_t:\ w\in C^{0,1}(\bom)
\to w_t\in C^{0,1}(\bom)$. By the a priori estimate (7.9), the
degree $\text{deg} (T_t, B_R, 0)$ is well defined, where $B_R$ is
the set of all positive function satisfying
$\|u\|_{C^{0,1}(\bom)}\le R$. When $t=0$, from (7.12) we have
obviously $w\equiv 1$. Namely $T_0$ has a unique fixed point
$w\equiv 1$. Hence the degree $\text{deg} (T_t, B_R, 0)=1$ for all
$t\in [0,1]$. This completes the proof. $\square$

\vskip10pt

\noo{\it Remark}. Theorem 1.2 extends to more general equations
(1.4) where
$$w=[\det D^2 u]^{\th-1}, \ \ \ 0<\th \le \frac 1n. $$

 \vskip30pt

\baselineskip=12pt
\parskip=1pt

\Refs\widestnumber\key{ABC}

\vskip10pt

\item {[1]}   W. Blaschke,
              Vorlesungen \'uber Differential geometrie,
              Berlin, 1923.

\item {[2]} L.A. Caffarelli,
             Interior a priori estimates for solutions of fully
             nonlinear equations,
             Ann. of Math. (2) 130 (1989),  189--213.

\item {[3]} L.A. Caffarelli,
             A localization property of viscosity solutions to
             the Monge-Amp\`ere equation and their strict convexity,
             Ann. Math., 131(1990), 129-134.

\item {[4]} L.A. Caffarelli,
             Interior $W^{2,p}$ estimates for solutions
             of Monge-Amp\`ere equations,
             Ann. Math., 131(1990), 135-150.

\item {[5]} L.A. Caffarelli,
             Boundary regularity of maps with convex potentials,
             Comm. Pure Appl. Math. 45 (1992), 1141--1151.

\item {[6]} L.A. Caffarelli,
             Boundary regularity of maps with convex potentials II.
             Ann. of Math. (2) 144 (1996),   453--496.

\item {[7]} L.A. Caffarelli and C.E. Guti\'errez,
             Properties of the solutions of the linearized
             Monge-Amp\`ere equations,
             Amer. J. Math., 119(1997), 423-465.

\item {[8]} L.A. Caffarelli, L. Nirenberg, and J. Spruck,
             The Dirichlet problem for nonlinear second order
             elliptic equations I. Monge-Amp\`ere equation,
             Comm. Pure Appl. Math., 37(1984), 369-402.

\item {[9]} E. Calabi,
             Improper affine hyperspheres of convex type and
             a generalization of a theorem by K. J\"orgens,
             Michigan Math. J. 5 1958 105--126.

\item {[10]} E. Calabi,
              Hypersurfaces with maximal affinely invariant area,
              Amer. J. Math. 104(1982), 91-126.

\item {[11]} S.Y. Cheng and S.T. Yau,
             Complete affine hypersurfaces, I.
             The completeness of affine metrics,
             Comm. Pure Appl. Math., 39(1986), 839-866.

\item {[12]} S.S. Chern,
             Affine minimal hypersurfaces,
             in {\it minimal submanifolds and geodesics},
             (Proc. Japan-United States Sem., Tokyo, 1977, 17-30.

\item {[13]} D. Gilbarg and N.S. Trudinger,
             Elliptic partial differential equations of second order,
             Springer-Verlag, New York, 1983.

\item {[14]} P. Guan, N.S. Trudinger, and X.-J. Wang,
             On the Dirichlet problem for degenerate
             Monge-Amp\`ere equations,
             Acta Math. 182 (1999), 87--104.

\item {[15]}  N.Ivochkina,
             A priori estimate of $\|u\|_{C^2(\bom)}$ of convex
             solutions of the Dirichlet problem for the
             Monge-Amp\`ere equation.
             Zap. Nauchn. Sem. Leningrad. Otdel. Mat. Inst. Steklov.
             (LOMI) 96 (1980), 69--79 (Russian). English
             translation in J. Soviet Math., 21(1983), 689-697.

\item {[16]}   N.V. Krylov,
              Nonlinear elliptic and parabolic equations of the
              second order, Reidel, Dordrecht-Boston, 1987.

\item {[17]}  K. Nomizu and T. Sasaki,
              Affine differential geometry,
              Cambridge University Press, 1994.

\item {[18]}    A.V. Pogorelov,
               The muitidimensional Minkowski problems,
               J. Wiley, New York, 1978.

\item {[19]}   M.V. Safonov,
              Classical solution of second-order nonlinear elliptic
              equations, Izv. Akad. Nauk SSSR Ser. Mat. 52 (1988),
              1272--1287 (Russian).  English translation in
              Math. USSR-Izv. 33 (1989), 597--612.

\item {[20]} N.S. Trudinger and X.-J. Wang,
              The Bernstein problem for affine maximal hypersurfaces,
              Invent. Math., 140 (2000), 399--422.

\item {[21]} N.S. Trudinger and X.-J. Wang,
              The affine Plateau problem,
              J. Amer. Math. Society, 18 (2005), 253-289.

\item {[22]}  X.-J. Wang,
              Regularity for Monge-Amp\`ere equation
              near the boundary,
              Analysis 16 (1996), 101--107.

\item {[23]}  X.-J. Wang,
              Affine maximal hypersurfaces,
              Proc. ICM  Vol.I\!I\!I, 2002, 221-231.

\endRefs

\enddocument

\end